\documentclass[a4paper,12pt]{amsart}
\pdfoutput=1
\usepackage{amsfonts, amsmath, amssymb, amstext}
\usepackage{tikz}
\usepackage[normalem]{ulem}
\usepackage{mathrsfs}
\usepackage{hyperref}
\usepackage{enumitem}

\allowdisplaybreaks[2]

\numberwithin{equation}{section}

\theoremstyle{plain}
\newtheorem{thm}{Theorem}[section]
\newtheorem{lem}[thm]{Lemma}
\newtheorem{prop}[thm]{Proposition}

\theoremstyle{definition}
\newtheorem{defn}[thm]{Definition}
\newtheorem{exam}[thm]{Example}

\theoremstyle{remark}
\newtheorem{rem}[thm]{Remark}

\newcommand{\norm}[1]{\left\|#1\right\|}
\newcommand{\cspan}{\overline{\mathrm{span}}}
\def\C{{\mathbb C}}
\def\N{{\mathbb N}}
\def\Q{{\mathbb Q}}

\def\Z{{\mathbb Z}}

\begin{document}

\title{Groupoid algebras as covariance algebras}

\author{Lisa Orloff Clark}
\address {Lisa Orloff Clark\\School of Mathematics and Statistics, Victoria University of Wellington, PO Box 600, Wellington 6140, New Zealand}
\email {lisa.clark@vuw.ac.nz}

\author{James Fletcher}
\address {James Fletcher\\School of Mathematics and Statistics, Victoria University of Wellington, PO Box 600, Wellington 6140, New Zealand}
\email {james.fletcher@vuw.ac.nz}

\date{\today}
\thanks{This research was supported by the Marsden grant VUW1514 from the
Royal Society of New Zealand.}



\begin{abstract}

Suppose $\mathcal{G}$ is a second-countable locally compact Hausdorff \'{e}tale groupoid, $G$ is a discrete group containing a unital subsemigroup $P$, and $c:\mathcal{G}\rightarrow G$ is a continuous cocycle. We derive conditions on the cocycle such that the reduced groupoid $C^*$-algebra $C_r^*(\mathcal{G})$ may be realised naturally as the covariance algebra of a product system over $P$ with coefficient algebra $C_r^*(c^{-1}(e))$. When $(G,P)$ is a quasi-lattice ordered group, we also derive conditions that allow $C_r^*(\mathcal{G})$ to be realised as the Cuntz--Nica--Pimsner algebra of a compactly aligned product system.
\end{abstract}

\maketitle

\section{Introduction}

Groupoid $C^*$-algebras generalise several important $C^*$-algebraic subclasses.  Many of the groupoid models that appear in the literature come with a natural homomorphism from the groupoid into a group, see for example \cite{MR2045419, MR1745529, MR1432596, MR2301938}.  In what follows, we use this kind of structure to realise a broad class of \'etale groupoid $C^*$-algebras as covariance algebras of natural product systems with coefficients in the groupoid $C^*$-subalgebra associated to the fibre of the identity.  

Fowler first introduced product systems of Hilbert bimodules in \cite{MR1907896}. Loosely speaking, a product system over a semigroup $P$ with coefficient algebra $A$ is a semigroup $\mathbf{X}=\bigsqcup_{p\in P} \mathbf{X}_p$, such that each $\mathbf{X}_p$ is a Hilbert $A$-bimodule, and the map $x\otimes_A y\mapsto xy$ extends to an isomorphism from $\mathbf{X}_p\otimes_A \mathbf{X}_q$ to $\mathbf{X}_{pq}$ for each $p,q\in P\setminus \{e\}$. Fowler focused on product systems over quasi-lattice ordered groups that satisfied a condition he called compact alignment. He then studied representations of such product systems satisfying a constraint called Nica covariance, as well as the associated universal $C^*$-algebra $\mathcal{NT}_\mathbf{X}$. Fowler also proposed a notion of Cuntz--Pimsner covariance for representations of product systems, and investigated the associated universal $C^*$-algebra $\mathcal{O}_\mathbf{X}$. Whilst Fowler's notion of Cuntz--Pimsner covariance coincides with the traditional notion of Cuntz--Pimsner covariance for a single Hilbert bimodule \cite{MR2102572, MR1426840}, in general the $C^*$-algebra $\mathcal{O}_\mathbf{X}$ fails to have a number of desirable properties that one might expect from a Cuntz--Pimsner like algebra. For example, Fowler's Cuntz--Pimsner algebra need not contain a faithful copy of the coefficient algebra $A$. Moreover, even when $\mathcal{O}_\mathbf{X}$ does contain a faithful copy of $A$, examples in the appendix of \cite{MR2069786} show that a representation of $\mathcal{O}_\mathbf{X}$ that is faithful on the copy of $A$ need not be faithful on the generalised fixed-point algebra associated to the canonical gauge coaction (this is a key step in establishing a gauge-invariant uniqueness theorem). Furthermore, in contrast with the Cuntz--Pimsner algebra associated to a single Hilbert bimodule, Fowler's Cuntz--Pimsner algebra $\mathcal{O}_\mathbf{X}$ need not be a quotient of $\mathcal{NT}_\mathbf{X}$.

In an attempt to overcome these issues, Sims and Yeend introduced a new covariance relation called Cuntz--Nica--Pimsner covariance \cite{MR2718947}. They showed that the universal $C^*$-algebra for representations satisfying this covariance relation, denoted by $\mathcal{NO}_{\mathbf{X}}$, coincides with Katsura's Cuntz--Pimsner algebra when $P=\N$ \cite[Proposition~5.3]{MR2718947}, and coincides with Fowler's $\mathcal{O}_{\mathbf{X}}$ whenever $A$ acts faithfully and compactly on each $\mathbf{X}_p$ and $P$ is directed (in the sense that each pair of elements in $P$ has a common upper bound) \cite[Proposition~5.1]{MR2718947}. Provided $\mathbf{X}$ satisfies a technical condition called $\tilde{\phi}$-injectivity, which is automatic for a wide class of examples, $\mathcal{NO}_\mathbf{X}$ is a quotient of $\mathcal{NT}_\mathbf{X}$, and the universal Cuntz--Nica--Pimsner representation of $\mathbf{X}$ is isometric \cite[Theorem~4.1]{MR2718947} (and so $\mathcal{NO}_\mathbf{X}$ contains a faithful copy of $A$). Furthermore, if
\begin{enumerate}[label=\upshape(\roman*)]
\item
the left action of $A$ on each fibre of $\mathbf{X}$ is faithful, or
\item
$P$ is directed and $\mathbf{X}$ is $\tilde{\phi}$-injective,
\end{enumerate}
then, subject to an amenability assumption, $\mathcal{NO}_\mathbf{X}$ has a gauge-invariant uniqueness theorem \cite[Corollary~4.11]{MR2837016}. To illustrate the utility of their construction, Sims and Yeend showed that finitely aligned higher-rank graph algebras \cite{MR2069786}, as well as Crisp and Laca's boundary quotients of Toeplitz algebras \cite{MR2274018}, have realisations as Cuntz--Nica--Pimsner algebras.

Whilst the Cuntz--Nica--Pimsner algebra introduced by Sims and Yeend certainly seems to be the correct Cuntz--Pimsner like algebra for compactly aligned product systems over quasi-lattice ordered groups $(G,P)$ satisfying (i) or (ii), issues still remain if we do not require these extra constraints. For example, in \cite[Example~3.16]{MR2718947} a non $\tilde{\phi}$-injective product system is exhibited for which the associated Cuntz--Nica--Pimsner algebra does not contain a faithful copy of the coefficient algebra. To overcome these issues, Sehnem has recently introduced the notion of strong covariance \cite{SEHNEM2018}. For an arbitrary product system, the associated covariance algebra always contains a faithful copy of the coefficient algebra, and any representation of the covariance algebra that is faithful on this copy of the coefficient algebra is faithful on the generalised fixed-point algebra associated to the canonical gauge coaction \cite[Theorem~3.10(C3)]{SEHNEM2018}. Furthermore, in the situation where either (i) or (ii) is satisfied, Sehnem's covariance algebra and Sims and Yeend's Cuntz--Nica--Pimsner algebra are naturally isomorphic \cite[Proposition~4.6]{SEHNEM2018}. The covariance algebra also has the advantage of being defined for a product system over an arbitrary unital semigroup that embeds in a group (rather than just a positive cone of a quasi-lattice ordered group), allowing for a much wider class of examples. Illustrating the utility of the construction, Sehnem shows that Li's semigroup $C^*$-algebras \cite{MR2900468}, as well as Exel's crossed products by interaction groups \cite{MR2422015}, can be realised as covariance algebras.

In this paper we develop conditions for a reduced groupoid $C^*$-algebra to have a natural realisation as the covariance algebra of a product system.  We also investigate when a groupoid $C^*$-algebra can be realised as the Cuntz--Nica--Pimsner algebra of a compactly aligned product system. The advantage of having these realisations is the complementary knowledge that can be gained from the two descriptions. We were motivated by the results of \cite{MR3624126}, where it is shown that the reduced groupoid $C^*$-algebra of an unperforated $\Z$-graded groupoid may be realised as the Cuntz--Pimsner algebra of a Hilbert bimodule whose coefficient algebra is the reduced groupoid $C^*$-algebra of the zero graded part of the groupoid.

The starting data for our construction is a second-countable locally compact Hausdorff \'{e}tale groupoid $\mathcal{G}$ equipped with a continuous group valued cocycle $c:\mathcal{G}\rightarrow G$. Supposing that the group $G$ contains a unital subsemigroup $P$, and defining
\[\mathbf{X}(\mathcal{G})_p:=\overline{C_c(c^{-1}(p))}\subseteq C_r^*(\mathcal{G})\] for each $p\in P$, we find necessary and sufficient conditions on the cocycle for \[\mathbf{X}(\mathcal{G}):=\bigsqcup_{p\in P}\mathbf{X}(\mathcal{G})_p\] to be a product system over $P$ with coefficient algebra $C_r^*(c^{-1}(e))$ (Proposition~\ref{bimodule isomorphisms}). It is then straightforward to check that the inclusion of $\mathbf{X}(\mathcal{G})$ in $C_r^*(\mathcal{G})$ is a representation, and we find necessary and sufficient conditions on the range map of $\mathcal{G}$ for this representation to be strongly covariant (Proposition~\ref{necessary and sufficient conditions for strong covariance}). Finally, we find necessary and sufficient conditions for the induced homomorphism to be surjective (Proposition~\ref{induced homomorphism is surjective - covariance algebra version}), and argue that if the group $G$ is amenable, then the induced homomorphism is injective (Theorem~\ref{induced homomorphism is injective - covariance algebra version}).

In Section~\ref{the QLOG situation}, we specialise to the situation where $(G,P)$ is a quasi-lattice ordered group. In Lemma~\ref{sufficient conditions for Nica covariance}, we develop necessary and sufficient conditions for $\mathbf{X}(\mathcal{G})$ to be compactly aligned and for the inclusion of $\mathbf{X}(\mathcal{G})$ in $C_r^*(\mathcal{G})$ to be Nica covariant. Combining this with the results from Section~\ref{realising as covariance algebras}, we have conditions for when $C_r^*(\mathcal{G})$ may be realised naturally as a Cuntz--Nica--Pimsner algebra (Theorem~\ref{main result for CNP algebras}). When $(G,P)=(\Z,\N)$ our conditions coincide with those of Rennie, Robertson, and Sims in \cite{MR3624126}, and so our result generalises \cite[Theorem~12]{MR3624126}.

Finally, in Section~\ref{examples} we present some examples of groupoids that satisfy our conditions. In particular, we look at the path and boundary-path groupoids associated to a topological higher-rank graph, as well as the groupoid associated to a directed semigroup action. In the future, we hope to apply our results to more exotic examples. For example, in \cite[Proposition~2.23]{Bourne2019}, using the results of  \cite{MR3624126},  it is shown that the reduced groupoid $C^*$-algebra associated to a one-dimensional Delone set may be realised as a Cuntz--Pimsner algebra --- the authors then ask if the reduced groupoid $C^*$-algebra associated to a $k$-dimensional Delone set may be realised as the Cuntz--Nica--Pimsner algebra of a product system over $\N^k$.

\section{Background and preliminaries}

\subsection{Hilbert bimodules}

Let $A$ be a $C^*$-algebra and $X$ be a right inner-product $A$-module with inner product $\langle \cdot, \cdot \rangle_A:X\times X\rightarrow A$.
The formula $\norm{x}_X:=\norm{\langle x,x\rangle_A}_A^{1/2}$ defines a norm on $X$ by \cite[Proposition~1.1]{lance}, and we say that $X$ is a Hilbert $A$-module if $X$ is complete with respect to this norm.
%
We write $\mathcal{L}_A(X)$ for the $C^*$-algebra of adjointable operators on a Hilbert $A$-module $X$.

For each $x,y\in X$ there is an adjointable operator $\Theta_{x,y}$ defined by $\Theta_{x,y}(z):=x\cdot \langle y, z\rangle_A$. We call such operators generalised rank-one operators. The closed subspace 
\[\mathcal{K}_A(X):=\cspan\{\Theta_{x,y}:x,y\in X\},\] elements of which we call generalised compact operators, forms an essential ideal of $\mathcal{L}_A(X)$.

A Hilbert $A$-bimodule (also known as a $C^*$-correspondence) is a Hilbert $A$-module $X$ equipped with a left action of $A$ by adjointable operators, i.e. there exists a homomorphism $\phi:A\rightarrow \mathcal{L}_A(X)$. To simplify notation, we will often write $a\cdot x$ for $\phi(a)(x)$. Since each $\phi(a)$ is by definition adjointable, and so $A$-linear, we have that $a\cdot (x\cdot b)=(a\cdot x)\cdot b$ for each $a,b\in A$ and $x\in X$.

A particularly simple (and important) example occurs when $X:=A$. Letting $A$ act on $X$ by left and right multiplication, and equipping $X$ with the $A$-valued inner product $\langle a,b\rangle_A:=a^*b$, gives a Hilbert $A$-bimodule, which we denote by ${}_A A_A$. The map $\Theta_{a,b}\mapsto ab^*$ extends to an isomorphism from $\mathcal{K}_A({}_A A_A)$ to $A$, whilst $\mathcal{L}_A({}_A A_A)$ is isomorphic to the multiplier algebra of $A$, which we denote by $\mathcal{M}(A)$.

Every Hilbert $A$-module $X$ is nondegenerate in the sense that the span of the set $\{x\cdot a:x\in X, a\in A\}$ is dense in $X$. In particular, the Hewitt--Cohen--Blanchard factorisation theorem \cite[Proposition~2.31]{MR1634408} says that for each $x\in X$, there exists a unique $x'\in X$ such that $x=x'\cdot \langle x', x'\rangle_A$. In general, a Hilbert $A$-bimodule need not be left nondegenerate in the sense that $X=\cspan\{a\cdot x:x\in X,a\in A\}$.

We can combine two Hilbert $A$-bimodules $X$ and $Y$ by taking their balanced tensor product. We let $X\odot Y$ denote the algebraic tensor product of $X$ and $Y$ as complex vector spaces, and write $X\odot_A Y$ for the quotient by 
\[\mathrm{span}\{x\cdot a \odot y- x\odot a\cdot y:x\in X, \, y\in Y, \, a\in A\}.\] Writing $x\odot_A y$ for the coset containing $x\odot y$, the formula 
\[\langle x\odot_A y, w\odot_A z\rangle_A:=\langle y, \langle x,w\rangle_A \cdot z\rangle_A\] determines a bounded $A$-valued sesquilinear form on $X\odot_A Y$. If we let $N$ be the subspace $\mathrm{span}\{n\in X\odot_A Y:\langle n,n\rangle_A=0\}$, then the formula  
\[\norm{z+N}:=\inf_{n\in N} \norm{\langle z+n,z+n\rangle_A}_A^{1/2}\] defines a norm on $(X\odot_A Y)/N$. The balanced tensor product of $X$ and $Y$, denoted by $X\otimes_A Y$, is then defined to be the completion of $(X\odot_A Y)/N$ with respect to this norm.

Another way to combine Hilbert bimodules is to take their direct sum. Given an indexing set $K$ and a collection $\{X_k:k\in K\}$ of Hilbert $A$-bimodules, we let $\bigoplus_{k\in K} X_k$ denote the space of sequences $(x_k)_{k\in K}$ such that $x_k\in X_k$ for each $k\in K$ and $\sum_{k\in K}\langle x_k,x_k\rangle_A$ converges in $A$. Proposition~1.1 of  \cite{lance} can be used to show that there exists an $A$-valued inner product on $\bigoplus_{k\in K} X_k$ such that 
\[\left\langle (x_k)_{k\in K}, (y_k)_{k\in K}\right\rangle_A=\sum_{k\in K}\langle x_k, y_k\rangle_A,\] and $\bigoplus_{k\in K} X_k$ is complete with respect to the induced norm. Letting $A$ act pointwise from the left and right gives $\bigoplus_{k\in K} X_k$ the structure of a Hilbert $A$-bimodule.



\subsection{Product systems and their representations}

\begin{defn}
Let $P$ be a semigroup with identity $e$, and $A$ a $C^*$-algebra. A product system over $P$ with coefficient algebra $A$ is a semigroup $\mathbf{X}=\bigsqcup_{p\in P}\mathbf{X}_p$ such that:
\begin{enumerate}[label=\upshape(\roman*)]
\item for each $p\in P$, $\mathbf{X}_p\subseteq \mathbf{X}$ is a Hilbert $A$-bimodule;
\item $\mathbf{X}_e$ is equal to the Hilbert $A$-bimodule ${}_A A_A$;
\item for each $p,q\in P\setminus \{e\}$, there exists a Hilbert $A$-bimodule isomorphism $M_{p,q}:\mathbf{X}_p\otimes_A \mathbf{X}_q\rightarrow \mathbf{X}_{pq}$ satisfying $M_{p,q}(x\otimes_A y)=xy$ for each $x\in \mathbf{X}_p$ and $y\in \mathbf{X}_q$; and
\item[(iv)] multiplication in $\mathbf{X}$ by elements of $\mathbf{X}_e=A$ implements the left and right actions of $A$ on each $\mathbf{X}_p$; that is $xa=x\cdot a$ and $ax=a\cdot x$ for each $p\in P$, $a\in A$, and $x\in \mathbf{X}_p$.
\end{enumerate}
\end{defn}

We write $\phi_p:A\rightarrow \mathcal{L}_A(\mathbf{X}_p)$ for the homomorphism that implements the left action of $A$ on $\mathbf{X}_p$, i.e. $\phi_p(a)(x)=a\cdot x=ax$ for each $p\in P$, $a\in A$, and $x\in \mathbf{X}_p$. Multiplication in $\mathbf{X}$ is associative since $\mathbf{X}$ is a semigroup. Hence, $\phi_{pq}(a)(xy)=(\phi_p(a)x)y$ for all $p,q\in P$, $a\in A$, $x\in \mathbf{X}_p$, and $y\in \mathbf{X}_q$. We write $M_{e,p}$ and $M_{p,e}$ for the maps that implement the left and right actions of $A$ on $\mathbf{X}_p$, i.e. $M_{e,p}(a\otimes_A x)=\phi_p(a)(x)$ and $M_{p,e}(x\otimes_A a)=x\cdot a$. Both $M_{e,p}$ and $M_{p,e}$ are inner-product preserving (and so injective), whilst $M_{p,e}$ is surjective by the Hewitt--Cohen--Blanchard factorisation theorem. Note: $M_{e,p}$ need not be surjective, since $\mathbf{X}_p$ is not necessarily nondegenerate as a left $A$-module. We also write $\langle \cdot,\cdot\rangle_A^p$ for the $A$-valued inner-product on $\mathbf{X}_p$.

For each $p\in P\setminus \{e\}$ and $q\in P$, we define a homomorphism $\iota_p^{pq}:\mathcal{L}_A\left(\mathbf{X}_p\right)\rightarrow \mathcal{L}_A\left(\mathbf{X}_{pq}\right)$ by
\[
\iota_p^{pq}(S):=M_{p,q}\circ (S\otimes_A \mathrm{id}_{\mathbf{X}_q})\circ M_{p,q}^{-1}
\]
for each $S\in \mathcal{L}_A\left(\mathbf{X}_p\right)$. Equivalently, $\iota_p^{pq}$ is characterised by the formula $\iota_p^{pq}(S)(xy)=(Sx)y$ for each $S\in \mathcal{L}_A\left(\mathbf{X}_p\right)$, $x\in \mathbf{X}_p$, $y\in \mathbf{X}_q$. For notational simplicity, we also define $\iota_e^q:=\phi_q$ as a map from $A\cong \mathcal{K}_A(\mathbf{X}_e)$ to $\mathcal{L}_A(\mathbf{X}_q)$. Given $p,r\in P$, we also define $\iota_p^r$ to be the zero map if $r\not\in pP$.

In order to associate $C^*$-algebras to product systems, we need the notion of a representation of a product system.

\begin{defn}
Let $\mathbf{X}$ be a product system over $P$ with coefficient algebra $A$. A representation of $\mathbf{X}$ in a $C^*$-algebra $B$ is a map $\psi:\mathbf{X}\rightarrow B$ satisfying the following relations:
\begin{enumerate}
\item[(T1)] each $\psi_p:=\psi|_{\mathbf{X}_p}$ is a linear map, and $\psi_e$ is a $C^*$-homomorphism;
\item[(T2)] $\psi_p(x)\psi_q(y)=\psi_{pq}(xy)$ for all $p,q\in P$ and $x\in \mathbf{X}_p$, $y\in \mathbf{X}_q$; and
\item[(T3)] $\psi_p(x)^*\psi_p(y)=\psi_e(\langle x,y\rangle_A^p)$ for all $p\in P$ and $x,y\in \mathbf{X}_p$.
\end{enumerate}
\end{defn}

Given a representation $\psi:\mathbf{X}\rightarrow B$, \cite[Proposition~8.11]{MR2135030} gives the existence of a homomorphism $\psi^{(p)}:\mathcal{K}_A\left(\mathbf{X}_p\right)\rightarrow B$ for each $p\in P$, such that $\psi^{(p)}\left(\Theta_{x,y}\right)=\psi_p(x)\psi_p(y)^*$ for all $x,y\in \mathbf{X}_p$.

The Toeplitz algebra of $\mathbf{X}$, which we denote by $\mathcal{T}_\mathbf{X}$, is the universal $C^*$-algebra for representations of $\mathbf{X}$. We denote the universal representation of $\mathbf{X}$ by $\tilde{t}$. Using the Hewitt--Cohen--Blanchard factorisation theorem, it is straightforward to show that
\[
\mathcal{T}_\mathbf{X}=\cspan\big\{\tilde{t}(x_{p_1})\tilde{t}(x_{p_2})^*\cdots \tilde{t}(x_{p_{n-1}})\tilde{t}(x_{p_n})^*:p_i\in P, x_{p_i}\in \mathbf{X}_{p_i}\big\}.
\]

Given a discrete group $G$, the universal property of the group $C^*$-algebra $C^*(G)$ induces a homomorphism $\delta_G:C^*(G)\rightarrow C^*(G)\otimes C^*(G)$ such that $\delta_G(i_G(g))=i_G(g)\otimes i_G(g)$ for each $g\in G$ (we use an unadorned $\otimes$ to denote the minimal tensor product of $C^*$-algebras). A (full) coaction of $G$ on a $C^*$-algebra $A$ is then an injective homomorphism $\delta:A\rightarrow A\otimes C^*(G)$, that satisfies the coaction identity
\[
(\delta\otimes \mathrm{id}_{C^*(G)})\circ \delta=(\mathrm{id}_A\otimes \delta_G)\circ \delta,
\]
and is nondegenerate in the sense that $A\otimes C^*(G)=\cspan\{\delta(A)(1_{\mathcal{M}(A)}\otimes C^*(G))\}$ (where $1_{\mathcal{M}(A)}$ is the identity of the multiplier algebra of $A$). A coaction of $G$ on $A$ provides a topological $G$-grading of $A$ (see for example \cite{Ruybook}), for which the spectral subspace at $g\in G$ is $A^g:=\{a\in A:\delta(a)=a\otimes i_G(g)\}$. We refer to the spectral subspace $A^e$ as the generalised fixed-point algebra.

If $\mathbf{X}$ is a product system over a semigroup $P$, and $P$ sits inside a group $G$, then the universal property of $\mathcal{T}_\mathbf{X}$ induces a coaction $\tilde{\delta}:\mathcal{T}_\mathbf{X}\rightarrow \mathcal{T}_\mathbf{X}\otimes C^*(G)$ such that
\[
\tilde{\delta}(\tilde{t}_p(x))=\tilde{t}_p(x)\otimes i_G(p) \quad\text{for each $p\in P$, $x\in \mathbf{X}_p$.}
\]
We call $\tilde{\delta}$ the generalised gauge coaction of $G$ on $\mathcal{T}_\mathbf{X}$. One can show that the spectral subspace $\mathcal{T}_\mathbf{X}^g$ at $g\in G$ is the closure of sums of the form
\[
\tilde{t}(x_{p_1})\tilde{t}(x_{p_2})^*\cdots \tilde{t}(x_{p_{n-1}})\tilde{t}(x_{p_n})^*
\]
where $p_1,\ldots,p_n\in P$ are such that $p_1p_2^{-1}\cdots p_{n-1}p_n^{-1}=g$, and $x_{p_i}\in \mathbf{X}_{p_i}$ for $i\in \{1,\ldots,n\}$.

\subsection{Strongly covariant representations}
We recall the construction of \cite{SEHNEM2018}.
Suppose that $\mathbf{X}$ is a product system over a semigroup $P$ with coefficient algebra $A$, and $P$ sits inside a group $G$. Given a finite set $F\subseteq G$, we define
\[
K_F:=\bigcap_{g\in F} gP,
\]
where $gP:=\{gh:h\in P\}$. For $p\in P$, we define an ideal $I_{p^{-1}(p\vee F)}\lhd A$ as follows: Firstly, for $g\in F$, let
\[
I_{p^{-1}K_{\{p,g\}}}:=
\begin{cases}
\bigcap_{r\in K_{\{p,g\}}}\mathrm{ker}(\phi_{p^{-1}r}) &\text{if $K_{\{p,g\}}\neq \emptyset$ and $p\not\in gP$}\\
A &\text{otherwise.}
\end{cases}
\]
We then set
\[
I_{p^{-1}(p\vee F)}:=\bigcap_{g\in F}I_{p^{-1}K_{\{p,g\}}}.
\]
We use these ideals to define two new Hilbert $A$-bimodules. Firstly, we let
\[
\mathbf{X}_F:=\bigoplus_{p\in P}\mathbf{X}_p\cdot I_{p^{-1}(p\vee F)}.
\]
In the proof of \cite[Proposition~3.5]{SEHNEM2018} it is shown that the left action of $A$ on $\mathbf{X}_F$ is always faithful.
We also have the following useful characterisation of elements of $\mathbf{X}_F$. The proof is almost exactly the same as that of \cite[Lemma~3.2]{MR2718947}.

\begin{lem}
\label{non extendable}
Let $G$ be a group and $P\subseteq G$ a unital subsemigroup. Suppose $\mathbf{X}$ is a product system over $P$ with coefficient algebra $A$.
Let $F$ be a finite subset of $G$ and suppose $x\in \mathbf{X}_p$ for some $p\in P$. Then $x\in \mathbf{X}_p\cdot I_{p^{-1}(p\vee F)}$ if and only if $xy=0$ whenever $g\in F$ is such that $pP\cap gP\neq \emptyset$, $p\not\in gP$, and $y\in \mathbf{X}_{p^{-1}r}$ for some $r\in pP\cap gP$.
\end{lem}

Secondly, writing $gF:=\{gh:h\in F\}$ for each $g\in G$, we define
\[
\mathbf{X}_F^+:=\bigoplus_{g\in G}\mathbf{X}_{gF}.
\]
The product system $\mathbf{X}$ has a natural representation on $\mathcal{L}_A(\mathbf{X}_F^+)$, which we now describe. For $p\in P$ and $x\in \mathbf{X}_p$, there exists an operator $t^F_p(x):\mathbf{X}_F^+\rightarrow \mathbf{X}_F^+$ such that
\[
((t^F_p(x)(y)_g)_q=
\begin{cases}
x(y_{p^{-1}g})_{p^{-1}q} &\text{if $q\in pP$}\\
0 &\text{otherwise}
\end{cases}
\]
for each $y\in \mathbf{X}_F^+$, $g\in G$, $q\in P$. For this operation to be well-defined, we need to know that $x(y_{p^{-1}g})_{p^{-1}q}\in \mathbf{X}_q\cdot I_{q^{-1}(q\vee gF)}$ whenever $q\in pP$. This follows from the fact that $(y_{p^{-1}g})_{p^{-1}q}\in \mathbf{X}_{p^{-1}q}\cdot I_{(p^{-1}q)^{-1}(p^{-1}q\vee p^{-1}gF)}$ and $I_{(p^{-1}q)^{-1}(p^{-1}q\vee p^{-1}gF)}=I_{q^{-1}(q\vee gF)}$ for $q\in pP$. We point out that the operator $t^F_p(x)$ maps the direct summand $\mathbf{X}_{gF}$ of $\mathbf{X}_F^+$ into the direct summand $\mathbf{X}_{pgF}$ for each $g\in F$. In particular, $t^F_p(x)$ maps the direct summand $\mathbf{X}_q\cdot I_{q^{-1}(q\vee gF)}$ of $\mathbf{X}_{gF}$ into the direct summand $\mathbf{X}_{pq}\cdot I_{(pq)^{-1}(pq\vee pgF)}$ of $\mathbf{X}_{pgF}$ for each $q\in P$.

To describe the adjoint of $t^F_p(x)$, we first define a map $\Theta_{x,q}^*:\mathbf{X}_{pq}\rightarrow \mathbf{X}_q$ for each $q\in P$. If $p=e$ (so that $x\in \mathbf{X}_e=A$), we let $\Theta_{x,q}^*:=\phi_q(x^*)$. On the other hand if $p\neq e$, then $\Theta_{x,q}^*$ is determined by the formula $\Theta_{x,q}^*(yw)=\langle x,y\rangle_A^p\cdot w$ for each $y\in \mathbf{X}_p$ and $w\in \mathbf{X}_q$. The adjoint of $t^F_p(x)$ is then given by
\[
(t^F_p(x)^*(y)_g)_q=\Theta_{x,q}^*((y_{pg})_{pq})
\]
for each $y\in \mathbf{X}_F^+$, $g\in G$, $q\in P$. For this operation to be well-defined, we need to know that $\Theta_x^q((y_{pg})_{pq})\in \mathbf{X}_q\cdot I_{q^{-1}(q\vee gF)}$. This follows from the fact that $(y_{pg})_{pq}\in \mathbf{X}_{pq}\cdot I_{(pq)^{-1}(pq\vee pgF)}$ and
$I_{(pq)^{-1}(pq\vee pgF)}=I_{q^{-1}(q\vee gF)}$. The adjoint $t^F_p(x)^*$ maps the direct summand $\mathbf{X}_{gF}$ of $\mathbf{X}_F^+$ into the direct summand $\mathbf{X}_{p^{-1}gF}$ for each $g\in F$. In particular, $t^F_p(x)^*$ maps the direct summand $\mathbf{X}_q\cdot I_{q^{-1}(q\vee gF)}$ of $\mathbf{X}_{gF}$ into the direct summand $\mathbf{X}_{p^{-1}q}\cdot I_{(p^{-1}q)^{-1}(p^{-1}q\vee p^{-1}gF)}$ of $\mathbf{X}_{p^{-1}gF}$ if $q\in pP$, and maps $\mathbf{X}_q\cdot I_{q^{-1}(q\vee gF)}$ to zero if $q\notin pP$.

Routine calculations show the collection of maps $t^F:=\{t^F_p\}_{p\in P}$ give a representation of $\mathbf{X}$ in $\mathcal{L}_A(\mathbf{X}_F^+)$, and so induces a homomorphism from $\mathcal{T}_\mathbf{X}$ to $\mathcal{L}_A(\mathbf{X}_F^+)$, which we denote by $t^F_*$.

For $g\in G$, we write $Q_g^F$ for the projection of $\mathbf{X}_F^+$ onto the direct summand $\mathbf{X}_{gF}$. It can be shown that
\[
t^F_p(x)Q_g^F=Q_{pg}^Ft^F_p(x) \quad \text{and} \quad t^F_p(x)^*Q_g^F=Q_{p^{-1}g}^Ft^F_p(x)^*
\]
for any $p\in P$ and $x\in \mathbf{X}_p$. We use these projections to define an ideal $J_e$ of the generalised fixed-point algebra $\mathcal{T}_\mathbf{X}^e$ by
\[
J_e:=\big\{b\in \mathcal{T}_\mathbf{X}^e:\lim_F\|b\|_F=0\big\},
\]
where $\|b\|_F:=\|Q_e^Ft^F_*(b)Q_e^F\|_{\mathcal{L}_A(\mathbf{X}_F^+)}$ and the limit is taken over the directed set consisting of all finite subsets of $G$. We are now ready to give the definition of strong covariance:

\begin{defn}[\cite{SEHNEM2018}, Definition~3.2]
We say that a representation $\psi:\mathbf{X}\rightarrow B$ is strongly covariant if the induced homomorphism $\psi_*:\mathcal{T}_\mathbf{X}\rightarrow B$ vanishes on $J_e$.
\end{defn}

The observant reader might be concerned that the strong covariance of a representation depends on the choice of the group that we are embedding $P$ into. Lemma~3.9 of \cite{SEHNEM2018} shows that there is nothing to worry about: if $G$ and $H$ are groups that contain $P$ as a subsemigroup, then a representation is strongly covariant with respect to $G$ if and only if it is strongly covariant with respect to $H$.

Finally, we are ready to give the definition of the covariance algebra of a product system and state some of its important properties.

\begin{thm}[\cite{SEHNEM2018}, Theorem~3.10]
Let $\mathbf{X}$ be a product system over a unital semigroup $P$ with coefficient algebra $A$. Suppose that $P$ is embeddable into a group $G$. Then there is a $C^*$-algebra $A\times _{\mathbf{X}}P$, called the covariance algebra of $\mathbf{X}$, and a strongly covariant representation $j_\mathbf{X}:\mathbf{X}\rightarrow A\times _{\mathbf{X}}P$ such that:
\begin{enumerate}[label=\upshape(\roman*)]
\item
$A\times _{\mathbf{X}}P$ is generated as a $C^*$-algebra by the image of $j_\mathbf{X}$;
\item
if $\psi:\mathbf{X}\rightarrow B$ is any strongly covariant representation, then there exists a unique homomorphism $\widehat{\psi}:A\times _{\mathbf{X}}P\rightarrow B$ such that $\widehat{\psi}\circ j_\mathbf{X}=\psi$;
\item
the homomorphism $j_{\mathbf{X}_e}:A\rightarrow A\times _{\mathbf{X}}P$ is faithful;
\item
there is a coaction $\delta$ of $G$ on $A\times _{\mathbf{X}}P$ such that $\delta(j_{\mathbf{X}_p}(x))=j_{\mathbf{X}_p}(x)\otimes i_G(p)$ for each $p\in P$ and $x\in\mathbf{X}_p$;
\item
a homomorphism from $A\times _{\mathbf{X}}P$ to a $C^*$-algebra $B$ is faithful on the generalised fixed-point algebra $(A\times _{\mathbf{X}}P)^\delta$  if and only if it is faithful on $j_{\mathbf{X}_e}(A)$.
\end{enumerate}
\end{thm}

\subsection{Quasi-lattice ordered groups and Cuntz--Nica--Pimsner algebras}
\label{QLOG background}

Imposing additional structure on the semigroup $P$ allows us to consider representations of $\mathbf{X}$ satisfying additional constraints.

Recall that a quasi-lattice ordered group $(G,P)$ consists of a group $G$ and a subsemigroup $P$ of $G$ such that $P\cap P^{-1}=\{e\}$, and, with respect to the partial order on $G$ induced by $p\leq q \Leftrightarrow p^{-1}q\in P $, any two elements $p,q\in G$ which have a common upper bound in $P$ have a least common upper bound in $P$. It is straightforward to show that if two elements in $G$ have a least common upper bound in $P$, then this least common upper bound is unique. Hence, if it exists, we write $p\vee q$ for the least common upper bound of $p,q\in G$. For $p,q\in G$, we write $p\vee q=\infty$ if $p$ and $q$ have no common upper bound in $P$, and $p\vee q<\infty$ otherwise. We say that $P$ is directed if $p\vee q<\infty$ for every $p,q\in P$.

Suppose that $(G,P)$ is a quasi-lattice ordered group and $\mathbf{X}$ is a product system over $P$ with coefficient algebra $A$. We say that $\mathbf{X}$ is compactly aligned if whenever $p,q\in P$ with $p\vee q<\infty$, we have that
\[
\iota_p^{p\vee q}(S)\iota_q^{p\vee q}(T)\in \mathcal{K}_A(\mathbf{X}_{p\vee q})
\quad \text{for all $S\in \mathcal{K}_A(\mathbf{X}_p)$, $T\in \mathcal{K}_A(\mathbf{X}_q)$.}
\]
It is important to note that this condition does not imply that either $\iota_p^{p\vee q}(S)$ or $\iota_q^{p\vee q}(T)$ is compact.

We say that a representation $\psi$ of a compactly aligned product system $\mathbf{X}$ is Nica covariant if, for any $p,q\in P$ and $S\in \mathcal{K}_A(\mathbf{X}_p)$, $T\in \mathcal{K}_A(\mathbf{X}_q)$, we have
\begin{align*}
\psi^{(p)}(S)\psi^{(q)}(T)=
\begin{cases}
\psi^{(p\vee q)}\left(\iota_p^{p\vee q}(S)\iota_q^{p\vee q}(T)\right) & \text{if $p\vee q<\infty$}\\
0 & \text{otherwise.}
\end{cases}
\end{align*}

The Nica--Toeplitz algebra of a compactly aligned product system $\mathbf{X}$, which we denote by $\mathcal{NT}_\mathbf{X}$, is the universal $C^*$-algebra for Nica covariant representations of $\mathbf{X}$. Denoting the universal Nica covariant representation of $\mathbf{X}$ by $i_\mathbf{X}$, it follows from the Hewitt--Cohen--Blanchard factorisation theorem that
\[
\mathcal{NT}_\mathbf{X}=\cspan\left\{i_\mathbf{X}(x)i_\mathbf{X}(y)^*:x,y\in \mathbf{X}\right\}.
\]

In addition to Nica covariance, we can ask that a representation of a compactly aligned product system satisfies an additional constraint called Cuntz--Pimsner covariance. Formulating this additional covariance relation requires some additional background material and notation, which we now present.

For $p\in P$, define an ideal $\overline{I}_p\lhd A$ by
\[
\overline{I}_p:=
\begin{cases}
\bigcap_{e<r\leq p}\mathrm{ker}(\phi_r) & \text{if $p\neq e$}\\
A & \text{if $p=e$.}
\end{cases}
\]
We use these ideals to define another Hilbert $A$-bimodule
\[
\widetilde{\mathbf{X}}_p:=\bigoplus_{e\leq r\leq p}\mathbf{X}_r\cdot \overline{I}_{r^{-1}p}.
\]
For $s\in P\setminus \{e\}$, there exists a homomorphism $\widetilde{\iota}_s^{\,p}:\mathcal{L}_A(\mathbf{X}_s)\rightarrow \mathcal{L}_A(\widetilde{\mathbf{X}}_p)$ such that
\[
\widetilde{\iota}_s^{\,p}(T)=\Bigg(\bigoplus_{s\leq r\leq p}\iota_s^r(T)|_{\mathbf{X}_r\cdot \overline{I}_{r^{-1}p}}\Bigg)\oplus\Bigg(\bigoplus_{s\not\leq r\leq p}0_{\mathbf{X}_r\cdot \overline{I}_{r^{-1}p}}\Bigg) \quad \text{for $T\in \mathcal{L}_A(\mathbf{X}_s)$}.
\]
We also define a homomorphism $\widetilde{\iota}_e^{\,p}:A\cong\mathcal{K}_A(\mathbf{X}_e)\rightarrow \mathcal{L}_A(\widetilde{\mathbf{X}}_p)$ by
\[
\widetilde{\iota}_e^{\,p}(a)=\bigoplus_{e\leq r\leq p}\phi_r(a).
\]
We say that $\mathbf{X}$ is $\tilde{\phi}$-injective if, for each $p\in P$, the homomorphism $\widetilde{\iota}_e^{\,p}$ is injective. This technical condition is often automatic --- as shown in \cite[Lemma~3.15]{MR2718947}, a compactly aligned product system $\mathbf{X}$ over a quasi-lattice ordered group $(G,P)$ is $\tilde{\phi}$-injective if $\phi_p$ is injective for each $p\in P$ or $(G,P)$ satisfies the following property:
\begin{equation}
\label{maximal element condition}
\text{\parbox{.85\textwidth}{
If $S\subseteq P$ is nonempty and there exists $q\in P$ such that $p\leq q$ for all $p\in S$, then there exists $p\in S$ such that $p\not\leq p'$ for all $p'\in S\setminus \{p\}$.
}}
\end{equation}

Before, we can discuss Cuntz--Pimsner covariance, we require one last definition. Given a quasi-lattice ordered group $(G,P)$, we say that a predicate statement $\mathcal{P}(p)$ (where $p\in P$) is true for large $p$ if, given any $r\in P$, there exists $q\geq r$, such that $\mathcal{P}(p)$ is true whenever $p\geq q$.

Finally, we are ready to present the definition of Cuntz--Pimsner covariance originally formulated by Sims and Yeend \cite[Definition~3.9]{MR2718947}. As in \cite{MR2718947} we give a definition only in the situation where $\mathbf{X}$ is $\tilde{\phi}$-injective. We say that a representation $\psi:\mathbf{X}\rightarrow B$ is Cuntz--Pimsner covariant if, for any finite set $F\subseteq P$ and any choice of compact operators $\left\{T_s\in \mathcal{K}_A\left(\mathbf{X}_s\right):s\in F\right\}$, we have that
\[
\sum_{s\in F} \widetilde{\iota}_s^{\,p} (T_p)=0 \text{ for large $p$} \quad \Rightarrow \quad \sum_{s\in F} \psi^{(s)}(T_s)=0.
\]

We say that a representation is Cuntz--Nica--Pimsner covariant if it is both Nica covariant and Cuntz--Pimsner covariant. The Cuntz--Nica--Pimsner algebra of a compactly aligned $\tilde{\phi}$-injective product system $\mathbf{X}$, which we denote by $\mathcal{NO}_\mathbf{X}$, is the universal $C^*$-algebra for Cuntz--Nica--Pimsner covariant representations of $\mathbf{X}$.

In general, the relationship between the notions of Nica covariance, Cuntz--Pimsner covariance, and strong covariance is quite subtle. Indeed, given a compactly aligned $\tilde{\phi}$-injective product system over a quasi-lattice ordered group, the associated Cuntz--Nica--Pimsner algebra and covariance algebra need not coincide. However, as shown in \cite[Proposition~4.6]{SEHNEM2018}, if $\mathbf{X}$ is a compactly aligned $\tilde{\phi}$-injective product system over a quasi-lattice ordered group $(G,P)$ with coefficient algebra $A$ and either $\phi_p$ is injective for each $p\in P$, or $P$ is directed, then $\mathcal{NO}_\mathbf{X}$ and $A\times_\mathbf{X} P$ are canonically isomorphic.

\subsection{Groupoids and their $C^*$-algebras}

A groupoid $\mathcal{G}$ is a small category in which every morphism is invertible. We define maps $r,s:\mathcal{G}\rightarrow \mathcal{G}$, called the range and source maps, by $r(\gamma):=\gamma \gamma^{-1}$ and $s(\gamma):=\gamma^{-1}\gamma$. We define the unit space of $\mathcal{G}$ as $\mathcal{G}^{(0)}:=r(\mathcal{G})=s(\mathcal{G})$. A pair $\gamma,\mu\in \mathcal{G}$ is composable if and only if $s(\gamma)=r(\mu)$, and in this case $\gamma\mu\in \mathcal{G}$. We write $\mathcal{G}^{(2)}:=\{(\gamma,\mu)\in \mathcal{G}\times \mathcal{G}:s(\gamma)=r(\mu)\}$ for the collection of all composable pairs in $\mathcal{G}$.

A topological groupoid is a groupoid equipped with a topology such that the inversion and composition maps are continuous (using the subspace topology on $\mathcal{G}^{(2)}$ inherited from the product topology on $\mathcal{G}\times \mathcal{G}$). For our purposes we will only consider the situation where the topology is second-countable, locally compact, and Hausdorff. We say that a topological groupoid is \'{e}tale if the range map (or equivalently the source map) is a local homeomorphism. In particular, the range and source maps in an \'{e}tale groupoid are open maps. An open bisection of $\mathcal{G}$ is an open subset $U\subseteq \mathcal{G}$ such that the restrictions $r|_U$, $s|_U$ are homeomorphisms onto open subsets of $\mathcal{G}^{(0)}$. It follows that $\mathcal{G}$ is \'{e}tale if and only if there exists a basis for the topology consisting of precompact (i.e. compact closure) open bisections.   Note that for precompact open bisections $U$ and $V$, the product 
\[UV = \{\gamma \mu : \gamma \in U, \mu \in V, s(\gamma)=r(\mu)\}\]
is also a precompact open bisection.

Given a second-countable locally compact Hausdorff \'{e}tale groupoid $\mathcal{G}$, we define the reduced groupoid $C^*$-algebra of $\mathcal{G}$ to be a completion of the space $C_c(\mathcal{G})$ of complex valued compactly supported functions on $\mathcal{G}$,  which we describe briefly below.  Notationally, if $U$ is a precompact open bisection, then we will view $C_c(U)$ as a subset of $C_c(\mathcal{G})$ by extending functions to take value 0 outside of $U$. 

Since $\mathcal{G}$ is \'{e}tale, $r^{-1}(u)$ and $s^{-1}(u)$ are discrete for any $u\in \mathcal{G}^{(0)}$. Thus, if $f\in C_c(\mathcal{G})$ and $u\in \mathcal{G}^{(0)}$, then $\mathrm{supp}(f)\cap r^{-1}(u)$ is finite. Consequently, we can define a convolution product and involution on $C_c(\mathcal{G})$ by
\[
(fg)(\gamma):=\sum_{\substack{\alpha \in \mathcal{G}:\\r(\alpha)=r(\gamma)}}f(\alpha)g(\alpha^{-1}\gamma)
\quad\text{and}\quad
f^*(\gamma):=\overline{f(\gamma^{-1})}
\]
for each $f,g\in C_c(\mathcal{G})$ and $\gamma\in \mathcal{G}$. These operations give $C_c(\mathcal{G})$ the structure of a $*$-algebra. Given $u\in \mathcal{G}^{(0)}$, there exists a representation $\pi_u: C_c(\mathcal{G})\rightarrow \mathcal{B}(\ell^2(s^{-1}(u)))$ such that
\[
\big(\pi_u(f)(\xi)\big)(\gamma):=\sum_{\substack{\alpha \in \mathcal{G}:\\r(\alpha)=r(\gamma)}}f(\alpha)\xi(\alpha^{-1}\gamma)
\]
for each $f\in C_c(\mathcal{G})$, $\xi\in \ell^2(s^{-1}(u))$, and $\gamma\in s^{-1}(u)$. We define the reduced groupoid $C^*$-algebra of $\mathcal{G}$, denoted by $C_r^*(\mathcal{G})$, to be the completion of $C_c(\mathcal{G})$ in the norm
\[
\|f\|_r:=\mathrm{sup}\big\{\|\pi_u(f)\|_{\mathcal{B}(\ell^2(s^{-1}(u)))}:u\in \mathcal{G}^{(0)}\big\}.
\]


\section{Constructing the product system}
\label{constructing the product system}

Suppose we have the following setup:
\begin{itemize}
\item
$\mathcal{G}$ is a second-countable locally compact Hausdorff \'{e}tale groupoid;
\item
$G$ is a group (equipped with the discrete topology) with unit $e$, and $P\subseteq G$ is a unital subsemigroup;
\item
$c:\mathcal{G}\rightarrow G$ is a continuous cocycle.
\end{itemize}
For each $g\in G$, we write $\mathcal{G}_g:=c^{-1}(g)$. Since $c$ is a cocycle, we have that $\mathcal{G}^{(0)}\subseteq \mathcal{G}_e$. As $c$ is continuous and $G$ has the discrete topology, each $\mathcal{G}_g$ is a clopen second-countable locally compact Hausdorff subspace of $\mathcal{G}$. Since $\mathcal{G}_e$ is also closed under multiplication and taking inverses, it forms an \'{e}tale subgroupoid of $\mathcal{G}$.

The goal of this paper is to find conditions on the groupoid $\mathcal{G}$ and the cocycle $c$ such that the reduced groupoid $C^*$-algebra $C_r^*(\mathcal{G})$ may be realised as the covariance algebra of a product system over the semigroup $P$ with coefficient algebra $C_r^*(\mathcal{G}_e)$. The situation when $G=\Z$ was investigated by Rennie, Robertson, and Sims in \cite{MR3624126}, and our results generalise theirs. 
By \cite[Theorem~3.10(C3)]{SEHNEM2018}, a necessary condition is that $C_r^*(\mathcal{G})$ contains a faithful copy of $C_r^*(\mathcal{G}_e)$.

\begin{lem}[\cite{MR3624126}, Lemma~3]
\label{faithful copy of coefficient algebra}
Let $A_e$ denote the completion of $C_c(\mathcal{G}_e)=\{f\in C_c(\mathcal{G}):\mathrm{supp}(f)\subseteq \mathcal{G}_e\}$ in $C_r^*(\mathcal{G})$. Then there exists an isomorphism $I_e:C_r^*(\mathcal{G}_e)\rightarrow A_e$ that extends the identity map on $C_c(\mathcal{G}_e)$.
\end{lem}

We will frequently use the homomorphism $I_e$ to identify $C_r^*(\mathcal{G}_e)$ with its image in $C_r^*(\mathcal{G})$. We define the fibres of our product system to be certain closed subspaces of $C_r^*(\mathcal{G})$. Before we give the definition of the fibres (and show that each has the structure of Hilbert $C_r^*(\mathcal{G}_e)$-bimodule), we state the following standard lemma.  The proof is straightforward.

\begin{lem}
\label{support of product}
For each $f,g\in C_c(\mathcal{G})$, $\mathrm{supp}(fg)\subseteq \mathrm{supp}(f)\mathrm{supp}(g)$.
\end{lem}

\begin{prop}
\label{defining the bimodules}
For each $p\in P\setminus \{e\}$, let $\mathbf{X}(\mathcal{G})_p$ denote the completion of $C_c(\mathcal{G}_p)=\{f\in C_c(\mathcal{G}):\mathrm{supp}(f)\subseteq \mathcal{G}_p\}$ in $C_r^*(\mathcal{G})$. Then
\begin{enumerate}[label=\upshape(\roman*)]
 \item\label{it1:definingbimodule}  $\mathbf{X}(\mathcal{G})_p \cap \mathbf{X}(\mathcal{G})_q = \{0\}$ for $p \neq q$; and 
 \item\label{it2:definingbimodule}
$\mathbf{X}(\mathcal{G})_p$ is a Hilbert $C_r^*(\mathcal{G}_e)$-module with right action and inner product given by
\begin{equation}
\label{right action and inner product}
\xi\cdot b=\xi I_e(b) \quad \text{and} \quad \langle \xi, \eta \rangle^p_{C_r^*(\mathcal{G}_e)}=I_e^{-1}(\xi^*\eta)
\quad \text{for $\xi,\eta\in \mathbf{X}(\mathcal{G})_p$ and $b\in C_r^*(\mathcal{G}_e)$.}
\end{equation}
The norm on $\mathbf{X}(\mathcal{G})_p$ induced by the inner product $\langle \cdot, \cdot \rangle^p_{C_r^*(\mathcal{G}_e)}$ agrees with the norm on $C_r^*(\mathcal{G})$. Furthermore, there exists a $*$-homomorphism $\phi_p:C_r^*(\mathcal{G}_e)\rightarrow \mathcal{L}_{C_r^*(\mathcal{G}_e)}(\mathbf{X}(\mathcal{G})_p)$ such that
\begin{equation}
\label{left action}
\phi_p(b)(\xi)=I_e(b)\xi  \quad \text{for $\xi\in \mathbf{X}(\mathcal{G})_p$ and $b\in C_r^*(\mathcal{G}_e)$,}
\end{equation}
giving $\mathbf{X}(\mathcal{G})_p$ the structure of a Hilbert $C_r^*(\mathcal{G}_e)$-bimodule.
\end{enumerate}
\begin{proof}
To prove \ref{it1:definingbimodule}, we use the continuous, injective, linear map 
\[j:C_r^*(\mathcal{G}) \to C_0(\mathcal{G}),\]
which is defined using the left-regular representation (see  \cite[Proposition~3.3.3]{Aidan}, for example).    
We also use that $j$ restricts to the identity map on $C_c(\mathcal{G})$.  
Fix $a \in  \mathbf{X}(\mathcal{G})_p \cap \mathbf{X}(\mathcal{G})_q$ with $p\neq q$.  Then there exist
sequences $\{f_n\} \subseteq C_c(\mathcal{G}_p)$ and 
$\{g_n\} \subseteq C_c(\mathcal{G}_q)$ that converge to $a$ in $C_r^*(\mathcal{G})$.  By continuity of $j$, both sequences converge to the function $j(a)$ (with respect to the uniform norm) in $C_0(\mathcal{G})$.  By way of contradiction, suppose that $a$ is nonzero.  Then the injectivity and linearity of $j$ implies that $j(a)$ is a nonzero function.  So there exists $\gamma \in \mathcal{G}$ such that $j(a)(\gamma) \neq 0$.  Thus, $f_n(\gamma)$ and $g_n(\gamma)$ are eventually nonzero and hence 
$\gamma \in \mathcal{G}_p \cap \mathcal{G}_q$, which is a contradiction.

For \ref{it2:definingbimodule}, it follows from Lemma~\ref{support of product} that if $p,q\in G$ and $f,g\in C_c(\mathcal{G})$ with $\mathrm{supp}(f)\subseteq \mathcal{G}_p$ and $\mathrm{supp}(g)\subseteq \mathcal{G}_q$, then $fg\in C_c(\mathcal{G})$ and $\mathrm{supp}(fg)\subseteq \mathcal{G}_p\mathcal{G}_q \subseteq \mathcal{G}_{pq}$ (because $c$ is a cocycle). Hence,
\[
C_c(\mathcal{G}_p)C_c(\mathcal{G}_q)\subseteq C_c(\mathcal{G}_{pq}) \quad \text{for each $p,q\in G$.}
\]
In particular
\[
C_c(\mathcal{G}_e)C_c(\mathcal{G}_p)C_c(\mathcal{G}_e)\subseteq C_c(\mathcal{G}_p) \quad \text{for each $p\in P$.}
\]
By continuity, it follows that if $a,b\in C_r^*(\mathcal{G}_e)$ and $\xi\in \mathbf{X}(\mathcal{G})_p$ , then  $I_e(a)\xi I_e(b)\in \mathbf{X}(\mathcal{G})_p$.

It is also routine to show that $C_c(\mathcal{G}_p)^*=C_c(\mathcal{G}_{p^{-1}})$ for each $p\in G$. This follows from the definition of the involution on $C_c(\mathcal{G})$ and the fact that $c$ is a cocycle: $f^*(\gamma)=\overline{f(\gamma^{-1})}$ and $c(\gamma^{-1})=c(\gamma)^{-1}$ for each $f\in C_c(\mathcal{G})$ and each $\gamma\in \mathcal{G}$. Thus,
\[
C_c(\mathcal{G}_p)^*C_c(\mathcal{G}_p)=C_c(\mathcal{G}_{-p})C_c(\mathcal{G}_p)\subseteq C_c(\mathcal{G}_e)
\quad \text{for each $p\in P$.}
\]
By continuity, it follows that if $\xi,\eta\in \mathbf{X}(\mathcal{G})_p$ , then  $\xi^*\eta\in I_e(C_r^*(\mathcal{G}_e))$.

Since $I_e$ is injective, it follows from the previous calculations and \cite[Lemma~3.2(1)]{MR1426840} that $\mathbf{X}(\mathcal{G})_p$ is a Hilbert $C_r^*(\mathcal{G}_e)$-module with right action and inner product given by \eqref{right action and inner product}. Moreover, \cite[Lemma~3.2(2)]{MR1426840} implies that for each $b\in C_r^*(\mathcal{G}_e)$, the map $\phi_p(b)$ defined by \eqref{left action} is adjointable. Since $I_e$ is a $*$-homomorphism, it follows that $\phi_p$ is also a $*$-homomorphism.
\end{proof}
\end{prop}

We want to show that the collection of Hilbert $C_r^*(\mathcal{G}_e)$-bimodules $\{\mathbf{X}(\mathcal{G})_p:p\in P\setminus \{e\}\}$ given by Proposition~\ref{defining the bimodules} gives a product system over $P$ with coefficient algebra $C_r^*(\mathcal{G}_e)$. To do this we will show that the multiplication of $\mathbf{X}(\mathcal{G})_p$ and $\mathbf{X}(\mathcal{G})_q$ in $C_r^*(\mathcal{G})$ gives a Hilbert $C_r^*(\mathcal{G}_e)$-bimodule isomorphism with $\mathbf{X}(\mathcal{G})_{pq}$ (for $p,q\in P\setminus \{e\}$). As we shall see, this map is always inner-product preserving (and hence injective), but need not be surjective in general. As such, we place an additional constraint on the cocycle $c$:
\begin{itemize}
\item
if $\gamma\in \mathcal{G}_{pq}$ for some $p,q\in P$, then there exist composable $\gamma'\in \mathcal{G}_{p}$ and $\gamma''\in \mathcal{G}_{q}$ such that $\gamma=\gamma'\gamma''$.
\end{itemize}
We call a cocycle satisfying this additional property unperforated. Note: if $P\subseteq c(r^{-1}(u))$ for each $u\in \mathcal{G}^{(0)}$ (this is almost saying that $c$ is strongly surjective), then $c$ is automatically unperforated. (To see this observe that if $\gamma\in \mathcal{G}_{pq}$, then $p\in c(r^{-1}(r(\gamma)))$. If we choose $\gamma'\in r^{-1}(r(\gamma))$ such that $c(\gamma')=p$ and let $\gamma'':=\gamma'^{-1}\gamma$, then we get what we need.)

In order to show that the multiplication map from $\mathbf{X}(\mathcal{G})_p\times \mathbf{X}(\mathcal{G})_q$ to $\mathbf{X}(\mathcal{G})_{pq}$ is surjective, we need the following preliminary result.  It is a generalisation of 
\cite[Lemma~7]{MR3624126}.

\begin{lem}
\label{factorisation result}
Suppose $c$ is an unperforated cocycle. Then for each $p,q\in P$, the space $\mathrm{span}\{gh:g\in C_c(\mathcal{G}_p),  h\in C_c(\mathcal{G}_q)\}$ is dense in $C_c(\mathcal{G}_{pq})$ in both the uniform norm and the bimodule norm from Proposition~\ref{defining the bimodules}.
\begin{proof}
Using an argument similar to  
\cite[Lemma~3.10]{MR2419901}, 
we see that for each $r \in P$,  $C_c(\mathcal{G}_{r})$ is equal to the span of functions $f \in C_c(\mathcal{G}_{r})$ such that the support of $f$ is contained in a precompact open bisection.  So to prove the lemma it suffices to show that if $U\subseteq \mathcal{G}_{pq}$ is a a fixed precompact open bisection, then
\[S:=\mathrm{span}\{gh:g\in C_c(\mathcal{G}_p),  h\in C_c(\mathcal{G}_q)  \text{ with }  \operatorname{supp}(gh) \subseteq U\}\] 
is dense in 
\[T:=\{f \in C_c(\mathcal{G}_{pq}): \operatorname{supp}(f) \subseteq U\}.\]
In $T$, the uniform norm and the bimodule norm (which is just the restriction of the norm on $C_r^*(\mathcal{G})$) coincide by \cite[Corollary~3.3.4]{Aidan}.  Thus, it is enough to show that $S$ is dense in $T$ in the uniform norm.
We do so using the Stone--Weierstrass theorem:  for distinct $x,y\in \mathcal{G}_{pq} \cap U$ we find functions $f_p\in C_c(\mathcal{G}_p)$ and $f_q\in C_c(\mathcal{G}_q)$  such that $f_pf_q \in S$,  $(f_pf_q)(x)=1$, and $(f_pf_q)(y)=0$ (this shows that $S$ separates points in $T$ and vanishes nowhere).

Since $c$ is unperforated, we can choose $x_p\in \mathcal{G}_p$ and $x_q\in \mathcal{G}_q$ such that $x=x_p x_q$.  Then since $\mathcal{G}$ has a basis of precompact open bisections, we can find precompact open bisections $U_q\subseteq \mathcal{G}_q$ containing $x_q$ and $U_p\subseteq \mathcal{G}_p$ containing $x_p$.  Since multiplication in $\mathcal{G}$ is continuous, we can arrange it so that $U_pU_q \subseteq U$ by taking intersections.  
Further, if $s(x)\neq s(y)$ then using the Hausdorff property in $\mathcal{G}^{(0)}$,  we can arrange it so that $s(y)\not \in s(U_q)$ by taking an intersection again.  

We claim that $y\not\in U_pU_q$.  First, suppose that $s(x)=s(y)$. Since $x\neq y$ and $s|_{U_pU_q}$ is injective, we see that $y\not\in U_pU_q$.  Otherwise, if $s(x)\neq s(y)$, then $s(y)\not\in s(U_q)$ by our choice of $U_q$. Since $s(U_p U_q)\subseteq s(U_q)$, we see that $y\not\in U_pU_q$ proving the claim.

Since $\mathcal{G}$ is locally compact and Hausdorff, it is completely regular, and we can separate closed sets from points in their complements with continuous functions (in fact $\mathcal{G}$ is normal because it is second-countable, and so we can separate disjoint closed sets with continuous functions). Since $x_p$ and $x_q$ are contained in the open sets $U_p$ and $U_q$ respectively, we can find continuous functions $f_p, f_q:\mathcal{G}\rightarrow [0,1]$ such $f_p(x_p)=1$, $f_q(x_q)=1$ and ${f_p}|_{\mathcal{G}\setminus U_p}\equiv 0$, ${f_q}|_{\mathcal{G}\setminus U_q}\equiv 0$. Since $U_p$ and $U_q$ are precompact, $f_p\in C_c(U_p)\subseteq C_c(\mathcal{G}_p)$ and $f_q\in C_c(U_q)\subseteq C_c(\mathcal{G}_q)$.  We also have 
$\operatorname{supp}(f_pf_q) \subseteq U$ by Lemma~\ref{support of product}. It remains to show that $(f_p f_q)(x)=1$ and $(f_p f_q)(y)=0$. We have
\[
(f_p f_q)(x)=(f_p f_q)(x_p x_q)=\sum_{\substack{\alpha \in \mathcal{G}:\\r(\alpha)=r(x_p x_q)}} f_p(\alpha) f_q(\alpha^{-1}x_px_q).
\]
Clearly, $r(x_p)=r(x_p x_q)$, and $f_p(x_p) f_q(x_p^{-1}x_px_q)=f_p(x_p) f_q(x_q)=1$. If $\alpha \in \mathcal{G}$ with $r(\alpha)=r(x_p x_q)$, then $r(\alpha)=r(x_p)$. If $\alpha \not\in U_p$ then $f_p(\alpha)=0$ since $\mathrm{supp}(f_p)\subseteq U_p$. On the other hand, if $\alpha\in U_p$, then $\alpha=x_p$ because $r|_{U_p}$ is injective. Thus, $(f_p f_q)(x)=1$. By Lemma~\ref{support of product}, $\mathrm{supp}(f_pf_q)\subseteq \mathrm{supp}(f_p)\mathrm{supp}(f_q)\subseteq U_p U_q$ , and so $(f_p f_q)(y)=0$ because $y\not\in U_pU_q$.
\end{proof}
\end{lem}

\begin{prop}
\label{bimodule isomorphisms}
For each $p,q\in P\setminus \{e\}$, there exists a $C_r^*(\mathcal{G}_e)$-linear inner-product preserving map $M_{p,q}:\mathbf{X}(\mathcal{G})_p\otimes_{C_r^*(\mathcal{G}_e)} \mathbf{X}(\mathcal{G})_q\rightarrow \mathbf{X}(\mathcal{G})_{pq}$ such that
\[
M_{p,q}(\xi\otimes_{C_r^*(\mathcal{G}_e)} \eta)=\xi\eta \quad \text{for $\xi \in \mathbf{X}(\mathcal{G})_p$ and $\eta \in \mathbf{X}(\mathcal{G})_q$}
\]
(where the multiplication on the right-hand side is occurring in $C_r^*(\mathcal{G})$). The maps $\{M_{p,q}:p,q\in P\setminus\{e\}\}$ are surjective if and only if $c$ is an unperforated cocycle.
\begin{proof}
We have already seen that $C_c(\mathcal{G}_p)C_c(\mathcal{G}_q)\subseteq C_c(\mathcal{G}_{pq})$ for each $p,q\in P$. Next observe that if $f,h\in C_c(\mathcal{G}_p)$ and $g,k\in C_c(\mathcal{G}_q)$, then
\begin{align*}
\langle f\otimes_{C_r^*(\mathcal{G}_e)} g, h\otimes_{C_r^*(\mathcal{G}_e)} k \rangle
=\langle g, \langle f, h \rangle_{C_r^*(\mathcal{G}_e)} \cdot k \rangle_{C_r^*(\mathcal{G}_e)}
=I_e^{-1}((fg)^*hk)
=\langle fg, hk \rangle^{pq}_{C_r^*(\mathcal{G}_e)}.
\end{align*}
It follows by linearity and continuity, that there exists an inner-product preserving map $M_{p,q}:\mathbf{X}(\mathcal{G})_p\otimes_{C_r^*(\mathcal{G}_e)} \mathbf{X}(\mathcal{G})_q\rightarrow \mathbf{X}(\mathcal{G})_{pq}$ such that $M_{p,q}(\xi\otimes_{C_r^*(\mathcal{G}_e)} \eta)=\xi\eta$ for each $\xi \in \mathbf{X}(\mathcal{G})_p$ and $\eta \in \mathbf{X}(\mathcal{G})_q$. Clearly, $M_{p,q}$ is left and right $C_r^*(\mathcal{G}_e)$-linear.

It remains to show that $M_{p,q}$ is surjective for each $p,q\in P\setminus \{e\}$ if and only if $c$ is unperforated. Firstly, suppose that $c$ is unperforated. By continuity, to show that $M_{p,q}$ is surjective it suffices to show that $C_c(\mathcal{G}_{pq})$ is contained in the range of $M_{p,q}$. This follows from the fact that $M_{p,q}$ is linear and isometric, and $\mathrm{span}\{gh:g\in C_c(\mathcal{G}_p), h\in C_c(\mathcal{G}_q)\}$ is dense in $C_c(\mathcal{G}_{pq})$ by Lemma~\ref{factorisation result}.

Conversely, suppose that $M_{p,q}$ is surjective for each $p,q\in P\setminus \{e\}$. Let $p,q\in P$ and suppose that $\gamma\in \mathcal{G}_{pq}$. If $p=e$, then $\gamma=r(\gamma)\gamma\in \mathcal{G}_p \mathcal{G}_q$. Similarly, if $q=e$, then $\gamma=\gamma s(\gamma)\in \mathcal{G}_p \mathcal{G}_q$. Hence, we need only worry about the situation when $p,q\neq e$. Choose a precompact open bisection $U\subseteq \mathcal{G}_{pq}$ containing $\gamma$, and $f\in C_c(U)\subseteq C_c(\mathcal{G}_{pq})\subseteq \mathbf{X}(\mathcal{G})_{pq}$ such that $f(\gamma)=1$. Since $M_{p,q}$ is surjective, we have that $f=\lim_{i\rightarrow \infty}\sum_{j_i=1}^{n_i}g_{j_i}h_{j_i}$ for some choice of $g_{j_i}\in C_c(\mathcal{G}_p)$ and $h_{j_i}\in C_c(\mathcal{G}_q)$. Since $f(\gamma)=1$, there must exist some $j_i$ such that
\[
0\neq (g_{j_i}h_{j_i})(\gamma)=\sum_{\substack{\alpha\in \mathcal{G}:\\r(\alpha)=r(\gamma)}}g_{j_i}(\alpha)h_{j_i}(\alpha^{-1}\gamma).
\]
Hence, there exists $\alpha\in \mathcal{G}$ with $r(\alpha)=r(\gamma)$ such that $\alpha\in \mathrm{supp}(g_{j_i})\subseteq \mathcal{G}_p$ and $\alpha^{-1}\gamma\in \mathrm{supp}(h_{j_i})\subseteq \mathcal{G}_q$. Thus, $\gamma=\alpha(\alpha^{-1}\gamma)\in \mathcal{G}_p\mathcal{G}_q$, and we conclude that $c$ is unperforated.
\end{proof}
\end{prop}

Proposition~\ref{bimodule isomorphisms} shows that if we let $\mathbf{X}(\mathcal{G})_e:={}_{C_r^*(\mathcal{G}_e)}C_r^*(\mathcal{G}_e)_{C_r^*(\mathcal{G}_e)}$ (which we will identify with $A_e\subseteq C_r^*(\mathcal{G})$ via $I_e$), then $\mathbf{X}(\mathcal{G}):=\bigsqcup_{p\in P}\mathbf{X}(\mathcal{G})_p$ equipped with the multiplication inherited from $C_r^*(\mathcal{G})$ is a product system over $P$ with coefficient algebra $C_r^*(\mathcal{G}_e)$.

\section{Groupoid $C^*$-algebras as covariance algebras}
\label{realising as covariance algebras}

Let $\mathbf{X}(\mathcal{G})$ be the product system over $P$ with coefficient algebra $C_r^*(\mathcal{G}_e)$ constructed in Section~\ref{constructing the product system}. For simplicity, we identify $C_r^*(\mathcal{G}_e)$ with $I_e(C_r^*(\mathcal{G}_e))\subseteq C_r^*(\mathcal{G})$. Our ultimate goal is to determine when the inclusion map $I:\mathbf{X}(\mathcal{G})\rightarrow C_r^*(\mathcal{G})$ induces an isomorphism from $C_r^*(\mathcal{G}_e)\times_{\mathbf{X}(\mathcal{G})}P$ to $C_r^*(\mathcal{G})$.

We begin by finding necessary and sufficient conditions for the inclusion of $\mathbf{X}(\mathcal{G})$ in $C_r^*(\mathcal{G})$ to be a strongly covariant representation. Routine calculations show that $I$ is always a representation, and so we just need to worry about strong covariance. We will shortly show that $I$ is strongly covariant  if and only if $\mathcal{G}$ satisfies the following condition:
for an odd integer $n \geq 1$ and $p,p_1\ldots,p_n\in P$,
\begin{equation}
\label{necessary and sufficient condition equivalent}
r(\mathcal{G}_p)\cap r(\mathcal{G}_{p_1}\mathcal{G}_{p_2}^{-1}\cdots \mathcal{G}_{p_{n-1}}^{-1}\mathcal{G}_{p_n})\subseteq \overline{\bigcup_{m\in pP\cap p_1p_2^{-1}\cdots p_{n-1}^{-1}p_nP}r(\mathcal{G}_m)}.
\end{equation}

The closure operation in \eqref{necessary and sufficient condition equivalent} plays two key roles. Firstly, the range map is in general not a closed map (it is always an open map however). Hence, despite each $\mathcal{G}_m$ being closed, $r(\mathcal{G}_m)$ need not be closed. Secondly, the set $pP\cap p_1p_2^{-1}\cdots p_{n-1}^{-1}p_nP$ need not be finite, and so even if we knew that each $r(\mathcal{G}_m)$ was closed, there would be no guarantee that $\bigcup_{m\in pP\cap p_1p_2^{-1}\cdots p_{n-1}^{-1}p_nP}r(\mathcal{G}_m)$ was also closed. We point out that in certain situations the union in \eqref{necessary and sufficient condition equivalent} can be replaced by a finite union. For example, if $(G,P)$ is a quasi-lattice ordered group, then for $p,p_1\ldots,p_n\in P$, we have
\[
pP\cap p_1p_2^{-1}\cdots p_{n-1}^{-1}p_nP
=
\begin{cases}
(p\vee p_1p_2^{-1}\cdots p_{n-1}^{-1}p_n)P & \text{if $p\vee p_1p_2^{-1}\cdots p_{n-1}^{-1}p_n<\infty$}\\
\emptyset & \text{otherwise,}
\end{cases}
\]
and so it follows that
\[
\bigcup_{m\in pP\cap p_1p_2^{-1}\cdots p_{n-1}^{-1}p_nP}r(\mathcal{G}_m)
=
\begin{cases}
r(\mathcal{G}_{p\vee p_1p_2^{-1}\cdots p_{n-1}^{-1}p_n}) & \text{if $p\vee p_1p_2^{-1}\cdots p_{n-1}^{-1}p_n<\infty$}\\
\emptyset & \text{otherwise}
\end{cases}
\]
because $c$ is unperforated.

The following result shows that in situations where the closure operation is not required to get the containment in \eqref{necessary and sufficient condition equivalent} (for example, see Example~\ref{semigroup action groupoids}), we have less to check.
\begin{lem}
\label{without closure, things are simpler}
The following conditions are equivalent
\begin{equation}
\label{a simpler useful inclusion}
r(\mathcal{G}_p)\cap r(\mathcal{G}_q)\subseteq \bigcup_{m\in pP\cap qP}r(\mathcal{G}_m) \quad \text{for all $p,q\in P$},
\end{equation}
\begin{equation}
\begin{aligned}
\label{a useful inclusion}
r(\mathcal{G}_p)\cap r(\mathcal{G}_{p_1}\mathcal{G}_{p_2}^{-1}\cdots \mathcal{G}_{p_{n-1}}^{-1}\mathcal{G}_{p_n})\subseteq &\bigcup_{m\in pP\cap p_1p_2^{-1}\cdots p_{n-1}^{-1}p_nP}r(\mathcal{G}_m)\\
&\quad\text{for all $n\geq1$ odd and $p,p_1,\ldots,p_n\in P$}.
\end{aligned}
\end{equation}
\begin{proof}
Clearly, \eqref{a useful inclusion} implies \eqref {a simpler useful inclusion}. Suppose that $\mathcal{G}$ satisfies \eqref{a simpler useful inclusion}. To prove that \eqref{a useful inclusion} holds, we will use induction on $n$. If $n=1$, then \eqref{a useful inclusion} is just \eqref{a simpler useful inclusion} and so there is nothing to prove. Now suppose that \eqref{a useful inclusion} holds for some odd $n\geq 1$. Fix $p,p_1,\ldots p_n,p_{n+1},p_{n+2}\in P$, and suppose that $z\in r(\mathcal{G}_p)\cap r(\mathcal{G}_{p_1}\mathcal{G}_{p_2}^{-1}\cdots \mathcal{G}_{p_{n-1}}^{-1}\mathcal{G}_{p_n}\mathcal{G}_{p_{n+1}}^{-1}\mathcal{G}_{p_{n+2}})$. Hence, there exist $y\in \mathcal{G}_p$, $w\in \mathcal{G}_{p_1}\mathcal{G}_{p_2}^{-1}\cdots \mathcal{G}_{p_{n-1}}^{-1}\mathcal{G}_{p_n}$, $u\in \mathcal{G}_{p_{n+1}}$, and $v\in \mathcal{G}_{p_{n+2}}$ such that $z=r(y)=r(wu^{-1}v)$. Then $z=r(y)=r(w)\in r(\mathcal{G}_p)\cap r(\mathcal{G}_{p_1}\mathcal{G}_{p_2}^{-1}\cdots \mathcal{G}_{p_{n-1}}^{-1}\mathcal{G}_{p_n})$ and $r(u)=r(v)\in r(\mathcal{G}_{p_{n+1}})\cap r(\mathcal{G}_{p_{n+2}})$. The inductive hypothesis says that there exists $m\in pP\cap p_1p_2^{-1}\cdots p_{n-1}^{-1}p_nP$ and $\alpha\in \mathcal{G}_m$ such that $z=r(\alpha)$. Similarly, there exists $n\in p_{n+1}P\cap p_{n+2}P$ and $\beta\in \mathcal{G}_n$ such that $r(u)=r(\beta)$. Then $s(w)=r(w^{-1}\alpha)=r(u^{-1}\beta)\in r(\mathcal{G}_{(p_1p_2^{-1}\cdots p_{n-1}^{-1}p_n)^{-1}m})\cap r(\mathcal{G}_{p_{n+1}^{-1}n})$, and so \eqref{a simpler useful inclusion} tells us that there exists $r\in (p_1p_2^{-1}\cdots p_{n-1}^{-1}p_n)^{-1}mP\cap p_{n+1}^{-1}nP$ and $\gamma\in \mathcal{G}_r$ such that $s(w)=r(\gamma)$. Since $z=r(w\gamma)\in r(\mathcal{G}_{p_1p_2^{-1}\cdots p_{n-1}^{-1}p_nr})$ and $p_1p_2^{-1}\cdots p_{n-1}^{-1}p_nr\in pP\cap p_1p_2^{-1}\cdots p_{n-1}^{-1}p_np_{n+1}^{-1}p_{n+2}P$, we conclude that \eqref{a useful inclusion} holds for arbitrary odd $n$.
\end{proof}
\end{lem}

Before proving that \eqref{necessary and sufficient condition equivalent} is necessary and sufficient for the inclusion map to be strongly covariant, we need a couple of preliminary results. The first is just a simple restatement of condition \eqref{necessary and sufficient condition equivalent} that is easier to work with.

\begin{lem}
\label{restatement of necessary and sufficient condition}
The groupoid $\mathcal{G}$ satisfies \eqref{necessary and sufficient condition equivalent} if and only if for all odd $n\geq1$ and $p,p_1\ldots,p_n\in P$, 
\begin{equation}
\label{necessary and sufficient condition}
\mathrm{int}\left(\left(r(\mathcal{G}_p)\cap r(\mathcal{G}_{p_1}\mathcal{G}_{p_2}^{-1}\cdots \mathcal{G}_{p_{n-1}}^{-1}\mathcal{G}_{p_n})\right) \setminus \bigcup_{m\in pP\cap p_1p_2^{-1}\cdots p_{n-1}^{-1}p_nP}r(\mathcal{G}_m)\right)=\emptyset.
\end{equation}
\begin{proof}
We will prove both directions of the lemma by contraposition. Firstly, suppose that $\mathcal{G}$ does not satisfy \eqref{necessary and sufficient condition equivalent}. Because $r$ is an open map, it follows that for some choice of $p,p_1\ldots,p_n\in P$,
\[
\left(r(\mathcal{G}_p)\cap r(\mathcal{G}_{p_1}\mathcal{G}_{p_2}^{-1}\cdots \mathcal{G}_{p_{n-1}}^{-1}\mathcal{G}_{p_n})\right)\setminus \overline{\bigcup_{m\in pP\cap p_1p_2^{-1}\cdots p_{n-1}^{-1}p_nP}r(\mathcal{G}_m)}
\]
is a nonempty open subset of
\[
\left(r(\mathcal{G}_p)\cap r(\mathcal{G}_{p_1}\mathcal{G}_{p_2}^{-1}\cdots \mathcal{G}_{p_{n-1}}^{-1}\mathcal{G}_{p_n})\right) \setminus \bigcup_{m\in pP\cap p_1p_2^{-1}\cdots p_{n-1}^{-1}p_nP}r(\mathcal{G}_m).
\]
Thus, \eqref{necessary and sufficient condition} does not hold.

Now suppose that $\mathcal{G}$ does not satisfy \eqref{necessary and sufficient condition}. Choose $p,p_1\ldots,p_n\in P$ and a nonempty open set $U\subseteq \left(r(\mathcal{G}_p)\cap r(\mathcal{G}_{p_1}\mathcal{G}_{p_2}^{-1}\cdots \mathcal{G}_{p_{n-1}}^{-1}\mathcal{G}_{p_n})\right) \setminus \bigcup_{m\in pP\cap p_1p_2^{-1}\cdots p_{n-1}^{-1}p_nP}r(\mathcal{G}_m)$. Looking for a contradiction, suppose that
$
U\not\subseteq \left(r(\mathcal{G}_p)\cap r(\mathcal{G}_{p_1}\mathcal{G}_{p_2}^{-1}\cdots \mathcal{G}_{p_{n-1}}^{-1}\mathcal{G}_{p_n})\right) \setminus \overline{\bigcup_{m\in pP\cap p_1p_2^{-1}\cdots p_{n-1}^{-1}p_nP}r(\mathcal{G}_m)}.
$
Hence, we can choose $x\in U$ and a sequence $(x_n)_{n=1}^\infty\subseteq \bigcup_{m\in pP\cap p_1p_2^{-1}\cdots p_{n-1}^{-1}p_nP}r(\mathcal{G}_m)$ such that $\lim_{n\rightarrow \infty}x_n=x$. Since $U$ is open, there exists $N\in \mathbb{N}$ such that $x_n\in U$ for all $n\geq N$.  But this is impossible by the choice of $U$. Thus,
$
\emptyset\neq U\subseteq \left(r(\mathcal{G}_p)\cap r(\mathcal{G}_{p_1}\mathcal{G}_{p_2}^{-1}\cdots \mathcal{G}_{p_{n-1}}^{-1}\mathcal{G}_{p_n})\right) \setminus \overline{\bigcup_{m\in pP\cap p_1p_2^{-1}\cdots p_{n-1}^{-1}p_nP}r(\mathcal{G}_m)}
$,
and we conclude that
\[
r(\mathcal{G}_p)\cap r(\mathcal{G}_{p_1}\mathcal{G}_{p_2}^{-1}\cdots \mathcal{G}_{p_{n-1}}^{-1}\mathcal{G}_{p_n})
\not\subseteq
\overline{\bigcup_{m\in pP\cap p_1p_2^{-1}\cdots p_{n-1}^{-1}p_nP}r(\mathcal{G}_m)}.
\]
Thus, $\mathcal{G}$ does not satisfy \eqref{necessary and sufficient condition equivalent}.
\end{proof}
\end{lem}

The second preliminary lemma forms the key part of our argument that $I$ is strongly covariant when $\mathcal{G}$ satisfies \eqref{necessary and sufficient condition equivalent}.

\begin{lem}
\label{useful technical lemma}
For $n\geq2$ even, let $p_1,\ldots,p_n\in P$ be such that  and $p_1p_2^{-1}\cdots p_{n-1}p_n^{-1}=e$. For each $i\in \{1,\ldots,n\}$, let $\xi_{p_i}\in C_c(\mathcal{G}_{p_i})\subseteq \mathbf{X}(\mathcal{G})_{p_i}$. Suppose that $F\subseteq G$ is finite and $p_1 p_2^{-1}\cdots p_{n-2k-2}^{-1}p_{n-2k-1}\in F$ for each $k\in \big\{0,\ldots,\frac{n-2}{2}\big\}$. If $\mathcal{G}$ satisfies \eqref{necessary and sufficient condition equivalent}, then
\begin{equation}
\label{a useful equality}
Q_e^Ft^F_*\big(\tilde{t}(\xi_{p_1})\tilde{t}(\xi_{p_2})^*\cdots\tilde{t}(\xi_{p_{n-1}})\tilde{t}(\xi_{p_{n}})^*\big)Q_e^F
=
t^F_e\big(\xi_{p_1}\xi_{p_2}^*\cdots \xi_{p_{n-1}}\xi_{p_n}^*\big)Q_e^F.
\end{equation}
\begin{proof}
Recall that for any $p\in P$, $\xi\in \mathbf{X}(\mathcal{G})_p$, and $g\in G$, we have
\begin{equation}
\label{interaction between t)F and Q_g^F}
t^F_p(\xi)Q_g^F=Q_{pg}^Ft^F_p(\xi) \quad \text{and} \quad t^F_p(\xi)^*Q_g^F=Q_{p^{-1}g}^Ft^F_p(\xi)^*.
\end{equation}
Since $p_1p_2^{-1}\cdots p_{n-1}p_n^{-1}=e$, repeatedly applying the two relations in \eqref{interaction between t)F and Q_g^F}, we see that
\begin{align*}
&Q_e^F t^F_*\big(\tilde{t}(\xi_{p_1})\tilde{t}(\xi_{p_2})^*\cdots\tilde{t}(\xi_{p_{n-1}})\tilde{t}(\xi_{p_{n}})^*\big)\\
&\qquad=Q_e^Ft^F_{p_1}(\xi_{p_1})t^F_{p_2}(\xi_{p_2})^*\cdots t^F_{p_{n-1}}(\xi_{p_{n-1}})t^F_{p_n}(\xi_{p_n})^*\\
&\qquad=t^F_{p_1}(\xi_{p_1})t^F_{p_2}(\xi_{p_2})^*\cdots t^F_{p_{n-1}}(\xi_{p_{n-1}})t^F_{p_n}(\xi_{p_n})^*Q_{p_np_{n-1}^{-1}\cdots p_2p_1^{-1}}^F\\
&\qquad=t^F_{p_1}(\xi_{p_1})t^F_{p_2}(\xi_{p_2})^*\cdots t^F_{p_{n-1}}(\xi_{p_{n-1}})t^F_{p_n}(\xi_{p_n})^*Q_{e}^F.
\end{align*}
Hence, for \eqref{a useful equality} to hold, it suffices to verify that
\begin{equation}
\begin{aligned}
\label{to be shown}
t^F_{p_1}(\xi_{p_1})t^F_{p_2}(\xi_{p_2})^*\cdots t^F_{p_{n-1}}(\xi_{p_{n-1}})t^F_{p_n}(\xi_{p_n})^*Q_{e}^F
=t^F_e\big(\xi_{p_1}\xi_{p_2}^*\cdots \xi_{p_{n-1}}\xi_{p_n}^*\big)Q_{e}^F.
\end{aligned}
\end{equation}
Let $x\in \mathbf{X}(\mathcal{G})_F^+$, and fix some $g\in G$ and $p\in P$. Then
\[
\big((t^F_e\big(\xi_{p_1}\xi_{p_2}^*\cdots \xi_{p_{n-1}}\xi_{p_n}^*\big)Q_{e}^Fx)_g\big)_p
=\begin{cases}
\xi_{p_1}\xi_{p_2}^*\cdots \xi_{p_{n-1}}\xi_{p_n}^*\big(x_e)_p & \text{if $g=e$}\\
0 & \text{otherwise.}
\end{cases}
\]
On the other hand,
\begin{align*}
\big((t^F_{p_1}(\xi_{p_1})&t^F_{p_2}(\xi_{p_2})^*\cdots t^F_{p_{n-1}}(\xi_{p_{n-1}})t^F_{p_n}(\xi_{p_n})^*Q_{e}^Fx)_g\big)_p\\
&=\begin{cases}
\xi_{p_1}\xi_{p_2}^*\cdots \xi_{p_{n-1}}\xi_{p_n}^*\big(x_e)_p & \text{if $g=e$ and $p\in p_1 p_2^{-1}\cdots p_{n-2k-2}^{-1}p_{n-2k-1}P$}\\
& \text{\qquad for each $k\in \big\{0,\ldots,\frac{n-2}{2}\big\}$}\\
0 & \text{otherwise.}
\end{cases}
\end{align*}
Thus, to establish \eqref{to be shown} it suffices to show that if $p\notin p_1 p_2^{-1}\cdots p_{n-2k-2}^{-1}p_{n-2k-1}P$ for some $k\in \big\{0,\ldots,\frac{n-2}{2}\big\}$, and $x\in \mathbf{X}(\mathcal{G})_p\cdot I_{p^{-1}(p\vee F)}$, then $\xi_{p_1}\xi_{p_2}^*\cdots \xi_{p_{n-1}}\xi_{p_n}^*x=0$.

Looking for a contradiction, suppose $p\in P$ is such that $p\notin p_1 p_2^{-1}\cdots p_{n-2k-2}^{-1}p_{n-2k-1}P$ for some $k\in \big\{0,\ldots,\frac{n-2}{2}\big\}$, and $x\in \mathbf{X}(\mathcal{G})_p\cdot I_{p^{-1}(p\vee F)}$ is such that $\xi_{p_1}\xi_{p_2}^*\cdots \xi_{p_{n-1}}\xi_{p_n}^*x\neq0$. Hence, there exists some $f\in C_c(\mathcal{G}_p)\cdot I_{p^{-1}(p\vee F)}$ such that $\xi_{p_1}\xi_{p_2}^*\cdots \xi_{p_{n-1}}\xi_{p_n}^*f\neq0$. Thus, $\mathrm{supp}(\xi_{p_1}\xi_{p_2}^*\cdots \xi_{p_{n-1}}\xi_{p_n}^*f)$ is a nonempty open subset of  $\mathcal{G}_p$.

We claim that
\begin{equation}
\label{range of support contained in intersection}
r(\mathrm{supp}(\xi_{p_1}\xi_{p_2}^*\cdots \xi_{p_{n-1}}\xi_{p_n}^*f))
\subseteq r(\mathcal{G}_p)
\cap
r(\mathcal{G}_{p_1}\mathcal{G}_{p_2}^{-1}\cdots \mathcal{G}_{p_{n-2k-2}}^{-1}\mathcal{G}_{p_{n-2k-1}}).
\end{equation}
Clearly, $r(\mathrm{supp}(\xi_{p_1}\xi_{p_2}^*\cdots \xi_{p_{n-1}}\xi_{p_n}^*f))\subseteq r(\mathcal{G}_p)$, so we just need to verify that
\begin{equation}
\label{range of support}
r(\mathrm{supp}(\xi_{p_1}\xi_{p_2}^*\cdots \xi_{p_{n-1}}\xi_{p_n}^*f))\subseteq r(\mathcal{G}_{p_1}\mathcal{G}_{p_2}^{-1}\cdots \mathcal{G}_{p_{n-2k-2}}^{-1}\mathcal{G}_{p_{n-2k-1}}).
\end{equation}
To see this, observe that if $z\in \mathrm{supp}(\xi_{p_1}\xi_{p_2}^*\cdots \xi_{p_{n-1}}\xi_{p_n}^*f)$, then
\begin{align*}
0&\neq (\xi_{p_1}\xi_{p_2}^*\cdots \xi_{p_{n-1}}\xi_{p_n}^*f)(z)\\
&=\sum_{\substack{\alpha\in \mathcal{G}:\\r(\alpha)=r(z)}}
(\xi_{p_1}\xi_{p_2}^*\cdots \xi_{p_{n-2k-2}}^*\xi_{p_{n-2k-1}})(\alpha)
(\xi_{p_{n-2k}}^*\xi_{p_{n-2k+1}}\cdots\xi_{p_{n-1}}\xi_{p_{n}}^*f)(\alpha^{-1}z).
\end{align*}
Hence, there exists $\alpha \in \mathcal{G}$ with $r(\alpha)=r(z)$ such that $\alpha\in \mathrm{supp}(\xi_{p_1}\xi_{p_2}^*\cdots \xi_{p_{n-2k-2}}^*\xi_{p_{n-2k-1}})$. Lemma~\ref{support of product} tells us that $\alpha\in \mathcal{G}_{p_1} \mathcal{G}_{p_2}^{-1}\cdots \mathcal{G}_{p_{n-2k-2}}^{-1}\mathcal{G}_{p_{n-2k-1}}$, and so $r(z)=r(\alpha)\in r(\mathcal{G}_{p_1} \mathcal{G}_{p_2}^{-1}\cdots \mathcal{G}_{p_{n-2k-2}}^{-1}\mathcal{G}_{p_{n-2k-1}})$. Since $z\in \mathrm{supp}(\xi_{p_1}\xi_{p_2}^*\cdots \xi_{p_{n-1}}\xi_{p_n}^*f)$ was arbitrary, we conclude that \eqref{range of support} holds

Next, we show that
\begin{equation}
\label{disjoint sets}
r(\mathrm{supp}(\xi_{p_1}\xi_{p_2}^*\cdots \xi_{p_{n-1}}\xi_{p_n}^*f))\cap\left(\bigcup_{m\in pP\cap p_1p_2^{-1}\cdots p_{n-2k-2}^{-1}p_{n-2k-1}P}r(\mathcal{G}_m)\right)=\emptyset.
\end{equation}
Looking for a contradiction, suppose that the intersection in \eqref{disjoint sets} is nonempty. Hence, we can choose $m\in pP\cap(p_1p_2^{-1}\cdots p_{n-2k-2}^{-1}p_{n-2k-1})P$, $\alpha\in \mathrm{supp}(\xi_{p_1}\xi_{p_2}^*\cdots \xi_{p_{n-1}}\xi_{p_n}^*f)$, and $\beta\in \mathcal{G}_m$ such that $r(\alpha)=r(\beta)$. Let $B\subseteq \mathcal{G}_{p^{-1}m}$ be a precompact open bisection containing $\alpha^{-1}\beta$. Choose a continuous function $\tau:\mathcal{G}\rightarrow [0,1]$ such that $\tau(\alpha^{-1}\beta)=1$ and $\tau|_{\mathcal{G}\setminus B}\equiv 0$. Then $\tau\in C_c(B)\subseteq C_c(\mathcal{G}_{p^{-1}m})\subseteq \mathbf{X}(\mathcal{G})_{p^{-1}m}$. Recall that $f\in C_c(\mathcal{G}_p)\cdot I_{p^{-1}(p\vee F)}$, and so $\xi_{p_1}\xi_{p_2}^*\cdots \xi_{p_{n-1}}\xi_{p_n}^*f\in \mathbf{X}(\mathcal{G})_p\cdot I_{p^{-1}(p\vee F)}$. Now because $p_1p_2^{-1}\cdots p_{n-2k-2}^{-1}p_{n-2k-1}\in F$ (one of the hypotheses of the lemma) and $p\not \in p_1p_2^{-1}\cdots p_{n-2k-2}^{-1}p_{n-2k-1}P$, Lemma~\ref{non extendable} thus tells us that $\xi_{p_1}\xi_{p_2}^*\cdots \xi_{p_{n-1}}\xi_{p_n}^*f\tau=0$. But this is impossible, because
\[
(\xi_{p_1}\xi_{p_2}^*\cdots \xi_{p_{n-1}}\xi_{p_n}^*f\tau)(\beta)=(\xi_{p_1}\xi_{p_2}^*\cdots \xi_{p_{n-1}}\xi_{p_n}^*f)(\alpha)\neq 0
\]
by construction. Thus, we conclude that \eqref{disjoint sets} holds.

Combining \eqref{range of support contained in intersection} and \eqref{disjoint sets}, we see that
\begin{align*}
r(\mathrm{supp}&(\xi_{p_1}\xi_{p_2}^*\cdots \xi_{p_{n-1}}\xi_{p_n}^*f))\\
&\subseteq \left(r(\mathcal{G}_p)\cap r(\mathcal{G}_{p_1}\mathcal{G}_{p_2}^{-1}\cdots \mathcal{G}_{p_{n-2k-2}}^{-1}\mathcal{G}_{p_{n-2k-1}})\right) \setminus \bigcup_{m\in pP\cap p_1p_2^{-1}\cdots p_{n-2k-2}^{-1}p_{n-2k-1}P}r(\mathcal{G}_m).
\end{align*}
Since $r(\mathrm{supp}(\xi_{p_1}\xi_{p_2}^*\cdots \xi_{p_{n-1}}\xi_{p_n}^*f))$ is nonempty and open (recall $r$ is a local homeomorphism, and so an open map), we see that
\[
\mathrm{int}\left(\left(r(\mathcal{G}_p)\cap r(\mathcal{G}_{p_1}\mathcal{G}_{p_2}^{-1}\cdots \mathcal{G}_{p_{n-2k-2}}^{-1}\mathcal{G}_{p_{n-2k-1}})\right) \setminus \bigcup_{m\in pP\cap p_1p_2^{-1}\cdots p_{n-2k-2}^{-1}p_{n-2k-1}P}r(\mathcal{G}_m)\right)
\neq\emptyset.
\]
However, in light of Lemma~\ref{restatement of necessary and sufficient condition}, this is impossible because $\mathcal{G}$ satisfies \eqref{necessary and sufficient condition equivalent}.

Hence, we conclude that if $p\notin p_1 p_2^{-1}\cdots p_{n-2k-2}^{-1}p_{n-2k-1}P$ for some $k\in \big\{0,\ldots,\frac{n-2}{2}\big\}$, and $x\in \mathbf{X}(\mathcal{G})_p\cdot I_{p^{-1}(p\vee F)}$, then $\xi_{p_1}\xi_{p_2}^*\cdots \xi_{p_{n-1}}\xi_{p_n}^*x=0$. Thus, \eqref{to be shown} holds, which proves the lemma.
\end{proof}
\end{lem}

Finally, we are ready to show that \eqref{necessary and sufficient condition equivalent} is necessary and sufficient for the inclusion map to be strongly covariant.

\begin{prop}
\label{necessary and sufficient conditions for strong covariance}
The inclusion map $I:\mathbf{X}(\mathcal{G})\rightarrow C_r^*(\mathcal{G})$ is a strongly covariant representation if and only if \eqref{necessary and sufficient condition equivalent} holds.
\begin{proof}
Firstly, suppose that \eqref{necessary and sufficient condition equivalent} holds. Let $b\in J_e$, i.e. $b\in \mathcal{T}_{\mathbf{X}(\mathcal{G})}^e$ and $\lim_F\|b\|_F=0$ (where $F$ ranges over the directed set determined by all finite subsets of $G$). Recall that $\|b\|_F:=\|Q^F_et^F_*(b)Q^F_e\|_{\mathcal{L}_A(\mathbf{X}(\mathcal{G})_F^+)}$. Write $I_*:\mathcal{T}_{\mathbf{X}(\mathcal{G})}\rightarrow C_r^*(\mathcal{G})$ for the homomorphism induced by the representation $I$.  In order for $I$ to be strongly covariant, we must show that $I_*(b)=0$.

Let $\varepsilon>0$. Since $b\in \mathcal{T}_{\mathbf{X}(\mathcal{G})}^e$, we can choose $N\in \mathbb{N}$, $n_j\in \mathbb{N}$ for $j\in \{1,\ldots N\}$, $p_i^j\in P$ for each $j\in \{1,\ldots N\}$ and $i\in \{1,\ldots,n_j\}$ such that $p_1^j(p_2^j)^{-1}\cdots p_{n_j-1}^j(p_{n_j}^j)^{-1}=e$, and $\xi_{p_i^j}\in C_c(\mathcal{G}_{p_i^j})$ for each $j\in \{1,\ldots N\}$ and $i\in \{1,\ldots,n_j\}$ such that
\[
\Big\|b-\sum_{j=1}^N \tilde{t}\big(\xi_{p_1^j}\big)\tilde{t}\big(\xi_{p_2^j}\big)^*\cdots\tilde{t}\big(\xi_{p_{n_j-1}^j}\big)\tilde{t}\big(\xi_{p_{n_j}^j}\big)^*\Big\|_{\mathcal{T}_{\mathbf{X}(\mathcal{G})}}<\varepsilon.
\]
For simplicity, we write $b':=\sum_{j=1}^N \tilde{t}\big(\xi_{p_1^j}\big)\tilde{t}\big(\xi_{p_2^j}\big)^*\cdots\tilde{t}\big(\xi_{p_{n_j-1}^j}\big)\tilde{t}\big(\xi_{p_{n_j}^j}\big)^*$.

Next choose a finite set $F\subseteq G$ such that if $F'\subseteq G$ is finite and $F\subseteq F'$, then $\|b\|_{F'}<\varepsilon$. For any such $F'$, we have
\begin{equation}
\label{norm estimate}
\begin{aligned}
\|b'\|_{F'}
&=\|Q^{F'}_et^F_*(b'-b+b)Q^{F'}_e\|_{\mathcal{L}_A(\mathbf{X}(\mathcal{G})_F^+)}\\
&\leq \|Q^{F'}_et^{F'}_*(b'-b)Q^{F'}_e\|_{\mathcal{L}_A(\mathbf{X}(\mathcal{G})_F^+)}+\|Q^{F'}_et^{F'}_*(b)Q^{F'}_e\|_{\mathcal{L}_A(\mathbf{X}(\mathcal{G})_F^+)}\\
&\leq\|t^{F'}_*(b'-b)\|_{\mathcal{L}_A(\mathbf{X}(\mathcal{G})_F^+)}+\|Q^{F'}_et^{F'}_*(b)Q^{F'}_e\|_{\mathcal{L}_A(\mathbf{X}(\mathcal{G})_F^+)} \quad \text{since $Q^{F'}_e$ is a projection}\\
&\leq\|b'-b\|_{\mathcal{T}_{\mathbf{X}(\mathcal{G})}}+\|Q^{F'}_et^{F'}_*(b)Q^{F'}_e\|_{\mathcal{L}_A(\mathbf{X}(\mathcal{G})_F^+)} \quad \text{since $t^{F'}_*$ is a representation of $\mathcal{T}_{\mathbf{X}(\mathcal{G})}$}\\
&=\|b'-b\|_{\mathcal{T}_{\mathbf{X}(\mathcal{G})}}+\|b\|_{F'}\\
&<2\varepsilon.
\end{aligned}
\end{equation}
By Lemma~\ref{useful technical lemma} we can choose a finite set $F'\subseteq G$ with $F\subseteq F'$ such that
\begin{equation}
\label{what we need to show}
Q^{F'}_et^F_*(b')Q^{F'}_e=t^{F'}_e\Big(\sum_{j=1}^N \xi_{p_1^j}\xi_{p_2^j}^*\cdots\xi_{p_{n_j-1}^j}\xi_{p_{n_j}^j}^*\Big)Q^{F'}_e.
\end{equation}
Combining \eqref{norm estimate} and \eqref{what we need to show}, and using the fact that $C_r^*(\mathcal{G}_e)$ acts faithfully on $\mathbf{X}(\mathcal{G})_{F'}$ (see the proof of \cite[Proposition~3.5]{SEHNEM2018}) for the second equality, we have that
\begin{align*}
\|I_*(b')\|_{C_r^*(\mathcal{G})}
&=\Big\|\sum_{j=1}^N \xi_{p_1^j}\xi_{p_2^j}^*\cdots\xi_{p_{n_j-1}^j}\xi_{p_{n_j}^j}^*\Big\|_{C_r^*(\mathcal{G})}\\
&=\Big\|t^{F'}_e\Big(\sum_{j=1}^N \xi_{p_1^j}\xi_{p_2^j}^*\cdots\xi_{p_{n_j-1}^j}\xi_{p_{n_j}^j}^*\Big)Q^{F'}_e\Big\|_{\mathcal{L}_A(\mathbf{X}(\mathcal{G})_F^+)}\\
&=\|Q^{F'}_et^F_*(b')Q^{F'}_e\|_{\mathcal{L}_A(\mathbf{X}(\mathcal{G})_F^+)}\\
&=\|b'\|_{F'}\\
&<2\varepsilon.
\end{align*}
Thus,
\[
\|I_*(b)\|_{C_r^*(\mathcal{G})}\leq \|I_*(b-b')\|_{C_r^*(\mathcal{G})}+\|I_*(b')\|_{C_r^*(\mathcal{G})}\leq \|b-b'\|_{\mathcal{T}_{\mathbf{X}(\mathcal{G})}}+\|I_*(b')\|_{C_r^*(\mathcal{G})}<3\varepsilon.
\]
Since $\varepsilon>0$ was arbitrary, we conclude that $I_*(b)=0$ as required. Thus, $I$ is strongly covariant.

Now we prove the converse by contraposition. Suppose that \eqref{necessary and sufficient condition equivalent} does not hold. By Lemma~\ref{restatement of necessary and sufficient condition}, we can choose $n\geq 1$ odd, $p,p_1\ldots,p_n\in P$, and a nonempty open set $U\subseteq \left(r(\mathcal{G}_p)\cap r(\mathcal{G}_{p_1}\mathcal{G}_{p_2}^{-1}\cdots \mathcal{G}_{p_{n-1}}^{-1}\mathcal{G}_{p_n})\right) \setminus \bigcup_{m\in pP\cap p_1p_2^{-1}\cdots p_{n-1}^{-1}p_nP}r(\mathcal{G}_m)$. Choose $y\in \mathcal{G}_p$ and $z_i\in \mathcal{G}_{p_i}$ such that $r(y)=r(z_1 z_2^{-1}\cdots z_{n-1}^{-1} z_n)\in U$. Choose a precompact open bisection $V\subseteq r^{-1}(U)\cap \mathcal{G}_p$ that contains $y$, and construct $f\in C_c(V)\subseteq C_c(\mathcal{G}_p)$ such that $f(y)=1$. Also choose precompact open bisections $V_i\subseteq \mathcal{G}_{p_i}$ such that $z_i\in V_i$, and construct $f_i\in C_c(V_i)\subseteq C_c(\mathcal{G}_{p_i})$ such that $f_i(z_i)=1$. Let
\[
b:=\tilde{t}(f)\tilde{t}(f)^*\tilde{t}(f_1)\tilde{t}(f_2)^*\cdots \tilde{t}(f_{n-1})^*\tilde{t}(f_n)\tilde{t}(f_n)^*\tilde{t}(f_{n-1})\cdots\tilde{t}(f_2)\tilde{t}(f_1)^*.
\]

We claim that $b\in J_e$. Firstly, using \cite[Lemma~2.2]{SEHNEM2018}, it is easy to see that $b\in \mathcal{T}_{\mathbf{X}(\mathcal{G})}^e$. Thus, it remains to show that $\lim_F \|b\|_F=0$ (taking the limit over all finite subsets of $G$). Let $F\subseteq G$ be finite. Looking for a contradiction, suppose that $t^F_*(b)Q^F_e\neq0$. Thus, there exists $q\in pP\cap p_1p_2^{-1}\cdots p_{n-1}^{-1}p_nP$ and $\tau\in C_c(\mathcal{G}_q)$ such that
\[
ff^*f_1f_2^*\cdots f_{n-1}^*f_nf_n^*f_{n-1}\cdots f_2f_1^*\tau\neq 0.
\]
Using Lemma~\ref{support of product}, we see that
\[
r(\mathrm{supp}(ff^*f_1f_2^*\cdots f_{n-1}^*f_nf_n^*f_{n-1}\cdots f_2f_1^*\tau))\subseteq r(V)\cap r(\mathcal{G}_q),
\]
and so $r(V)\cap r(\mathcal{G}_q)$ is nonempty. On the other hand,
\[
r(V)\subseteq U\subseteq \left(r(\mathcal{G}_p)\cap r(\mathcal{G}_{p_1}\mathcal{G}_{p_2}^{-1}\cdots \mathcal{G}_{p_{n-1}}^{-1}\mathcal{G}_{p_n})\right) \setminus \bigcup_{m\in pP\cap p_1p_2^{-1}\cdots p_{n-1}^{-1}p_nP}r(\mathcal{G}_m)
\]
by our choice of $V$ and $U$. Hence we have our contradiction, and we conclude that $t^F_*(b)Q^F_e=0$ for any finite subset $F\subseteq G$. Thus, $\lim_F \|b\|_F=\lim_F\|Q^F_et^F_*(b)Q_e^F\|_{\mathcal{L}_A(\mathbf{X}(\mathcal{G})_F^+)}=0$, and so $b\in J_e$ as claimed.

Finally, observe that
\[
I_*(b)=ff^*f_1f_2^*\cdots f_{n-1}^*f_nf_n^*f_{n-1}\cdots f_2 f_1^*,
\]
which is nonzero because
\[
(ff^*f_1f_2^*\cdots f_{n-1}^*f_nf_n^*f_{n-1}\cdots f_2 f_1^*)(r(y))=1
\]
by construction. Thus, $I_*$ does not vanish on $J_e$, and so $I$ is not strongly covariant.
\end{proof}
\end{prop}

When the inclusion map $I:\mathbf{X}(\mathcal{G})\rightarrow C_r^*(\mathcal{G})$ is a strongly covariant representation, we get a homomorphism $\widehat{I}:C_r^*(\mathcal{G}_e)\times_{\mathbf{X}(\mathcal{G})}P\rightarrow C_r^*(\mathcal{G})$ such that $\widehat{I}(j_{\mathbf{X}(\mathcal{G})_p}(x))=x$ for each $p\in P$ and $x\in \mathbf{X}(\mathcal{G})_p$. We now need to determine when $\widehat{I}$ is an isomorphism. Firstly, we find necessary and sufficient conditions for $\widehat{I}$ to be surjective.

\begin{prop}
\label{induced homomorphism is surjective - covariance algebra version}
The homomorphism $\widehat{I}:C_r^*(\mathcal{G}_e)\times_{\mathbf{X}(\mathcal{G})}P\rightarrow C_r^*(\mathcal{G})$ is surjective if and only if $c^{-1}(P)$ generates $\mathcal{G}$ as a groupoid.
\begin{proof}
Firstly, suppose that $\widehat{I}$ is surjective. Since every element of $C_r^*(\mathcal{G}_e)\times_{\mathbf{X}(\mathcal{G})}P$ can be approximated by sums of elements of the form
\[
j_{\mathbf{X}(\mathcal{G})_{p_1}}(\xi_{p_1})j_{\mathbf{X}(\mathcal{G})_{p_2}}(\xi_{p_2})^*
\cdots j_{\mathbf{X}(\mathcal{G})_{p_{n-1}}}(\xi_{p_{n-1}})j_{\mathbf{X}(\mathcal{G})_{p_n}}(\xi_{p_n})^*
\]
where $p_i\in P$ and $\xi_{p_i}\in C_c(\mathcal{G}_{p_i})$, whilst $\widehat{I}(j_{\mathbf{X}(\mathcal{G})_p}(\xi))=\xi$ for each $p\in P$ and $\xi\in C_c(\mathcal{G}_p)$, we see that every element of $C_r^*(\mathcal{G})$ can be approximated by sums of elements of the form
\[
\xi_{p_1}\xi_{p_2}^*\cdots \xi_{p_{n-1}}\xi_{p_n}^*
\]
where $p_i\in P$ and $\xi_{p_i}\in C_c(\mathcal{G}_{p_i})$. Thus, if $\gamma\in \mathcal{G}$, then there exist $p_i\in P$ and $\xi_{p_i}\in C_c(\mathcal{G}_{p_i})$ such that
\[
(\xi_{p_1}\xi_{p_2}^*\cdots \xi_{p_{n-1}}\xi_{p_n}^*)(\gamma)\neq 0.
\]
By Lemma~\ref{support of product}, it follows that
\begin{align*}
\gamma
\in \mathrm{supp}(\xi_{p_1})\mathrm{supp}(\xi_{p_2}^*)\cdots \mathrm{supp}(\xi_{p_{n-1}})\mathrm{supp}(\xi_{p_n}^*)
\subseteq \mathcal{G}_{p_1}\mathcal{G}_{p_2}^{-1}\cdots \mathcal{G}_{p_{n-1}}\mathcal{G}_{p_n}^{-1}.
\end{align*}
Thus, $\gamma$ belongs to the groupoid generated by $c^{-1}(P)$. Since $\gamma\in \mathcal{G}$ was arbitrary, we conclude that $c^{-1}(P)$ generates $\mathcal{G}$ as a groupoid.

Conversely, suppose that $c^{-1}(P)$ generates $\mathcal{G}$ as a groupoid. Since 
\[
\widehat{I}\big(C_r^*(\mathcal{G}_e)\times_{\mathbf{X}(\mathcal{G})}P\big)=\cspan\{\xi_{p_1}\xi_{p_2}^*\cdots \xi_{p_{n-1}}\xi_{p_n}^*:p_i\in P, \xi_{p_i}\in C_c(\mathcal{G}_{p_i})\}
\] 
and $C_c(\mathcal{G})$ is dense in $C_r^*(\mathcal{G})$, in order to show that $\widehat{I}$ is surjective, we need only show that $\mathrm{span}\{\xi_{p_1}\xi_{p_2}^*\cdots \xi_{p_{n-1}}\xi_{p_n}^*:p_i\in P, \xi_{p_i}\in C_c(\mathcal{G}_{p_i})\}$ is dense in $C_c(\mathcal{G})$. Like the proof of Lemma~\ref{factorisation result}, it suffices to show uniform density for functions supported on a precompact open bisection $U\subseteq \mathcal{G}$.  We show that if $x,y\in U$ are distinct, then there exist $p_i\in P$ and $\xi_{p_i}\in C_c(\mathcal{G}_{p_i})$ such that $\operatorname{supp}(\xi_{p_1}\xi_{p_2}^*\cdots \xi_{p_{n-1}}\xi_{p_n}^*) \subseteq U$,  $(\xi_{p_1}\xi_{p_2}^*\cdots \xi_{p_{n-1}}\xi_{p_n}^*)(x)=1$ and $(\xi_{p_1}\xi_{p_2}^*\cdots \xi_{p_{n-1}}\xi_{p_n}^*)(y)=0$.  Density will then follow by applying the Stone--Weierstrass Theorem.  

Since $c^{-1}(P)$ generates $\mathcal{G}$ as a groupoid and $P$ is a semigroup, there exist $p_i\in P$ and $x_i\in \mathcal{G}_{p_i}$ such that $x=x_1x_2^{-1}\ldots x_{n-1}x_n^{-1}$. As in Lemma~\ref{factorisation result}, we can choose precompact open bisections $U_i\subseteq \mathcal{G}_{p_i}$ with $x_i\in U_i$ such that $U_1U_2^{-1}\cdots U_{n-1}U_{n}^{-1} \subseteq U$ and if $s(x)\neq s(y)$, then $s(y)\not\in s(U_n^{-1})$. Next choose continuous functions $\xi_{p_i}:\mathcal{G}_{p_i}\rightarrow [0,1]$ such that $\xi_{p_i}(x_i)=1$ and $\xi_{p_i}|_{\mathcal{G}_{p_i}\setminus U_i}\equiv 0$. Thus, $\xi_{p_i}\in C_c(U_i)\subseteq C_c(\mathcal{G}_{p_i})$. We have that $\operatorname{supp}(\xi_{p_1}\xi_{p_2}^*\cdots \xi_{p_{n-1}}\xi_{p_n}^*) \subseteq U$ by Lemma~\ref{support of product}.  One can then check that $(\xi_{p_1}\xi_{p_2}^*\cdots \xi_{p_{n-1}}\xi_{p_n}^*)(x)=1$ and $(\xi_{p_1}\xi_{p_2}^*\cdots \xi_{p_{n-1}}\xi_{p_n}^*)(y)=0$ as required.
\end{proof}
\end{prop}

With the additional assumption that the group $G$ is amenable, it is routine to verify that $\widehat{I}$ is injective. Recall that a coaction $\eta$ of a group $H$ on a $C^*$-algebra $C$ is said to be normal if $(\mathrm{id}_C\otimes \lambda_H)\circ \eta$ is injective (where $\lambda_H:C^*(H)\rightarrow C_r^*(H)$ is the left regular representation). Observe that if $H$ is amenable, then every coaction of $H$ is normal because $\lambda_H$ is an isomorphism.

\begin{thm}
\label{induced homomorphism is injective - covariance algebra version}
Suppose that $G$ is an amenable group. Then the inclusion of $\mathbf{X}(\mathcal{G})$ in $C_r^*(\mathcal{G})$ induces an isomorphism from $C_r^*(\mathcal{G}_e)\times_{\mathbf{X}(\mathcal{G})}P$ to $C_r^*(\mathcal{G})$ if and only if $\mathcal{G}$ satisfies \eqref{necessary and sufficient condition equivalent} and $c^{-1}(P)$ generates $\mathcal{G}$ as a groupoid.
\begin{proof}
It follows from Proposition~\ref{necessary and sufficient conditions for strong covariance} and Proposition~\ref{induced homomorphism is surjective - covariance algebra version} that the inclusion of $\mathbf{X}(\mathcal{G})$ in $C_r^*(\mathcal{G})$ induces a surjective homomorphism from $C_r^*(\mathcal{G}_e)\times_{\mathbf{X}(\mathcal{G})}P$ to $C_r^*(\mathcal{G})$ if and only if $\mathcal{G}$ satisfies \eqref{necessary and sufficient condition equivalent} and $c^{-1}(P)$ generates $\mathcal{G}$ as a groupoid.

Thus, to complete the proof, it remains to show that if $G$ is amenable, and $\mathcal{G}$ satisfies \eqref{necessary and sufficient condition equivalent} and $c^{-1}(P)$ generates $\mathcal{G}$ as a groupoid, then the induced homomorphism $\widehat{I}:C_r^*(\mathcal{G}_e)\times_{\mathbf{X}(\mathcal{G})}P\rightarrow C_r^*(\mathcal{G})$ is injective. To do this we will use the abstract uniqueness theorem given by \cite[Corollary~A.3]{MR3047633}. Firstly, observe that $\widehat{I}\circ j_{\mathbf{X}(\mathcal{G})_e}=I_e$, which is injective by Lemma~\ref{faithful copy of coefficient algebra}. Thus, by \cite[Theorem~3.10(C3)]{SEHNEM2018}, $\widehat{I}$ is faithful on the fixed-point algebra $(C_r^*(\mathcal{G}_e)\times_{\mathbf{X}(\mathcal{G})}P)^\delta$.  Since $G$ is amenable, \cite[Lemma~6.1]{2017arXiv171101052M} gives the existence of a coaction $\beta:C_r^*(\mathcal{G})\rightarrow C_r^*(\mathcal{G})\otimes C^*(G)$ such that $\beta(f)=f\otimes i_G(g)$ for each $g\in G$ and $f\in C_c(\mathcal{G}_g)$ (note: we are making use of the amenability here because \cite[Lemma~6.1]{2017arXiv171101052M} only gives a reduced coaction). It is then routine to check that $\beta\circ \widehat{I}$ and $(\widehat{I}\otimes \mathrm{id}_{C^*(G)})\circ \delta$ agree on the generators of  $C_r^*(\mathcal{G}_e)\times_{\mathbf{X}(\mathcal{G})}P$, and so are equal as homomorphisms. Finally, since $G$ is amenable, the coactions $\delta$ and $\beta$ are normal. Hence, \cite[Corollary~A.3]{MR3047633} tells us that $\widehat{I}$ is injective.
\end{proof}
\end{thm}

\section{Specialising to the situation where $(G,P)$ is quasi-lattice ordered}
\label{the QLOG situation}

We now consider the situation where $(G,P)$ is a quasi-lattice ordered group. It is then natural to ask when the product system $\mathbf{X}(\mathcal{G})$ is compactly aligned and whether $C_r^*(\mathcal{G})$ is isomorphic to the Cuntz--Nica--Pimsner algebra of $\mathbf{X}(\mathcal{G})$.

\begin{lem}
\label{sufficient conditions for Nica covariance}
If
\begin{equation}
\label{extendibility}
r(\mathcal{G}_p)\cap r(\mathcal{G}_q)=
\begin{cases}
r(\mathcal{G}_{p\vee q}) & \text{if $p\vee q<\infty$}\\
\emptyset & \text{if $p\vee q=\infty$}
\end{cases}
\end{equation}
for each $p,q\in P$, then $\mathbf{X}(\mathcal{G})$ is compactly aligned. If $\mathbf{X}(\mathcal{G})$ is compactly aligned and the inclusion of $\mathbf{X}(\mathcal{G})$ in $C_r^*(\mathcal{G})$ is Nica covariant, then $\mathcal{G}$ satisfies \eqref{extendibility}.
\begin{proof}
Observe that since $c$ is unperforated, if $p,q\in P$ and $p\leq q$, then $r(\mathcal{G}_q)\subseteq r(\mathcal{G}_p)$. Hence, we always have that $r(\mathcal{G}_{p\vee q})\subseteq r(\mathcal{G}_p)\cap r(\mathcal{G}_q)$ for any $p,q\in P$ with $p\vee q<\infty$.

We now show that if $\mathcal{G}$ satisfies \eqref{extendibility}, then $\mathbf{X}(\mathcal{G})$ is compactly aligned. Let $p,q\in P$ with $p\vee q<\infty$ and suppose that $f,g\in C_c(\mathcal{G}_p)$ and $h,k\in C_c(\mathcal{G}_q)$. Routine calculations show that $\mathrm{supp}(fg^*hk^*)\subseteq \mathcal{G}_e\cap s^{-1}(r(\mathcal{G}_p)\cap r(\mathcal{G}_q))=\mathcal{G}_e\cap s^{-1}(r(\mathcal{G}_{p\vee q}))$ (where the last equality comes from \eqref{extendibility}). We claim that $\mathrm{span}\{mn^*:m,n\in C_c(\mathcal{G}_{p\vee q})\}$ is dense in $C_c(\mathcal{G}_e\cap s^{-1}(r(\mathcal{G}_{p\vee q})))$.

Once again, as in the proof of Lemma~\ref{factorisation result}, it suffices to show uniform density for functions supported on a precompact open bisection $U$.    Since $\mathcal{G}_e\cap s^{-1}(r(\mathcal{G}_{p\vee q}))\cap U$ is open, it is locally compact, and so we can use the Stone--Weierstrass Theorem as before.  Thus it suffices to show that for distinct $x,y\in s^{-1}(r(\mathcal{G}_{p\vee q}))\cap \mathcal{G}_e \cap U$ there exist $m,n\in C_c(\mathcal{G}_{p\vee q})$ such that $\operatorname{supp}(mn^*) \subseteq U$,  $(mn^*)(x)=1$ and $(mn^*)(y)=0$.  

Since $x\in s^{-1}(r(\mathcal{G}_{p\vee q}))\cap \mathcal{G}_e$, there exists $\alpha \in \mathcal{G}_{p\vee q}$ with $r(\alpha)=s(x)$, and for every such $\alpha$, we have that $x=(x\alpha)\alpha^{-1}\in \mathcal{G}_{p\vee q}\mathcal{G}_{p\vee q}^{-1}$. Choose a precompact open bisection $V\subseteq \mathcal{G}_{p\vee q}$ containing $\alpha$ such that if $s(x)\neq s(y)$, then $s(y)\not\in s(V^{-1})$.
We also choose a precompact open bisection $W$ inside $\mathcal{G}_{p\vee q}$ containing $x\alpha$.  By taking intersections, we can arrange it so that $WV^{-1} \subseteq U$. Now choose continuous functions $m,n:\mathcal{G}\rightarrow \C$ such that $m(x\alpha)=1$, $m|_{\mathcal{G}\setminus W}\equiv 0$ and $n(\alpha)=1$, $n|_{\mathcal{G}\setminus V}\equiv 0$. Since $W$ and $V$ are precompact, $m\in C_c(W)\subseteq C_c(\mathcal{G}_{p\vee q})$ and $n\in C_c(V)\subseteq C_c(\mathcal{G}_{p\vee q})$. By Lemma~\ref{support of product}, $\operatorname{supp}(mn^*) \subseteq U$.  We now show that $(mn^*)(x)=1$ and $(mn^*)(y)=0$. Firstly,
\[
(mn^*)(x)=\sum_{\substack{\gamma\in \mathcal{G}:\\r(\gamma)=r(x)}}m(\gamma)n(x^{-1}\gamma).
\]
Clearly, $r(x\alpha)=r(x)$. If $\gamma\in \mathcal{G}$ with $r(\gamma)=r(x)$ and $\gamma\in \mathrm{supp}(m)\subseteq W$, then $\gamma=x\alpha$ since $x\alpha\in W$ by construction and $W$ is a bisection. If $\gamma=x\alpha$, then $x^{-1}\gamma=xx^{-1}\alpha=\alpha$. Thus, $(mn^*)(x)=m(x\alpha)n(\alpha)=1$. Since $\mathrm{supp}(mn^*)\subseteq \mathrm{supp}(m)\mathrm{supp}(n)^{-1}\subseteq WV^{-1}$, in order to show that $(mn^*)(y)=0$, it suffices to show that $y\not\in WV^{-1}$. Looking for a contradiction, suppose that $y\in WV^{-1}$. Then $s(x)\neq s(y)$ since the source map is injective on the bisection $WV^{-1}$ (recall that the collection of bisections is closed under multiplication and taking inverses) and $x\neq y$. Thus, by our choice of $V$, $s(y)\not \in s(V^{-1})$. But this is impossible since $s(WV^{-1})\subseteq s(V^{-1})$. Hence, $(mn^*)(y)=0$. By the Stone--Weierstrass Theorem, we conclude that $\mathrm{span}\{mn^*:m,n\in C_c(\mathcal{G}_{p\vee q})\}$ is dense in $C_c(s^{-1}(r(\mathcal{G}_{p\vee q}))\cap \mathcal{G}_e)$.

Hence, we may choose $m_{j_i},n_{j_i}\in C_c(\mathcal{G}_{p\vee q})$ such that $fg^*hk^*=\lim_{i\rightarrow \infty} \sum_{j_i=1}^{n_j} m_{j_i}n_{j_i}^*$. It is then straightforward to verify that
\[
\iota_p^{p\vee q}(\Theta_{f,g})\iota_q^{p\vee q}(\Theta_{h,k})=\phi_{p\vee q}(fg^*hk^*)
=\lim_{i\rightarrow \infty} \sum_{j_i=1}^{n_j} \Theta_{m_{j_i},n_{j_i}}
\in \mathcal{K}_{C_r^*(\mathcal{G}_e)}(\mathbf{X}(\mathcal{G})_{p\vee q}).
\]
By linearity and continuity, it follows that $\mathbf{X}(\mathcal{G})$ is compactly aligned.

To finish the proof, it remains to show that if $\mathbf{X}(\mathcal{G})$ is compactly aligned and the representation of $\mathbf{X}(\mathcal{G})$ in $C_r^*(\mathcal{G})$ given by the inclusion map is Nica covariant, then \eqref{extendibility} holds. Let $p,q\in P$ and suppose $\xi\in r(\mathcal{G}_p)\cap r(\mathcal{G}_q)$. Choose $y\in \mathcal{G}_p$ and $z\in \mathcal{G}_q$ such that $\xi=r(y)=r(z)$ (so $\xi=yy^{-1}=zz^{-1}$). Fix precompact open bisections $U\subseteq \mathcal{G}_p$ and $V\subseteq \mathcal{G}_q$ containing $y$ and $z$ respectively, and choose continuous functions $f,g:\mathcal{G}\rightarrow [0,1]$ such that $f(y)=1=g(z)$ and $f|_{\mathcal{G}\setminus U}\equiv 0$ and $g|_{\mathcal{G}\setminus V}\equiv 0$. Then $f\in C_c(U)\subseteq C_c(\mathcal{G}_p)$ and $g\in C_c(V)\subseteq C_c(\mathcal{G}_q)$. One can then show that
\begin{equation}
\label{simple calculation}
I^{(p)}(\Theta_{f,f})I^{(q)}(\Theta_{g,g})(\xi)=(ff^*gg^*)(\xi)=1.
\end{equation}
If $p\vee q=\infty$, then the Nica covariance of $I$ tells us that $I^{(p)}(\Theta_{f,f})I^{(q)}(\Theta_{g,g})=0$, which is impossible. Hence, the intersection $r(\mathcal{G}_p)\cap r(\mathcal{G}_q)$ is empty whenever $p\vee q=\infty$.
Thus, to establish that \eqref{extendibility} holds, it remains to show that if $p\vee q<\infty$, then $\xi\in r(\mathcal{G}_{p\vee q})$. Since $\mathbf{X}(\mathcal{G})$ is compactly aligned, we can write $\iota_p^{p\vee q}(\Theta_{f,f})\iota_q^{p\vee q}(\Theta_{g,g})=\lim_{i\rightarrow \infty} \sum_{j_i=1}^{n_i} \Theta_{\mu_{j_i},\nu_{j_i}}$ for some choice of $\mu_{j_i},\nu_{j_i}\in C_c(\mathcal{G}_{p\vee q})$. Using the Nica covariance of $I$, we have
\begin{align*}
1
=I^{(p)}(\Theta_{f,f})I^{(q)}(\Theta_{g,g})(\xi)
&=I^{(p\vee q)}\Big(\lim_{i\rightarrow \infty} \sum_{j_i=1}^{n_i} \Theta_{\mu_{j_i},\nu_{j_i}}\Big)(\xi)\\
&=\lim_{i\rightarrow \infty} \sum_{j_i=1}^{n_i} \mu_{j_i}\nu_{j_i}^*(\xi)\\
&=\lim_{i\rightarrow \infty} \sum_{\substack{1\leq j_i\leq n_i\\ \alpha\in r^{-1}(\xi)}}\mu_{j_i}(\alpha)\overline{\nu_{j_i}(\alpha)}.
\end{align*}
Thus, there exists $j_i$ such that $\emptyset \neq r^{-1}(\xi)\cap \mathrm{supp}(\mu_{j_i})$. Since $\mathrm{supp}(\mu_{j_i})\subseteq \mathcal{G}_{p\vee q}$, we conclude that $\xi\in r(\mathcal{G}_{p\vee q})$, and so \eqref{extendibility} holds.
\end{proof}
\end{lem}

The following example shows that \eqref{extendibility} is not necessary if we just require $\mathbf{X}(\mathcal{G})$ to be compactly aligned. The example also shows that in general the inclusion of $\mathbf{X}(\mathcal{G})$ in $C_r^*(\mathcal{G})$ need not be Cuntz--Pimsner covariant (however, since $\mathcal{G}$ does not satisfy \eqref{extendibility}, Lemma~\ref{sufficient conditions for Nica covariance} tells us that $I$ is also not Nica covariant).

\begin{exam}
\label{when things can go wrong}
Define $\alpha:\Q_+^*\rightarrow \mathrm{Aut}(\Q)$ by $\alpha_h(g):=hg$. One can show that for any $(x,y), (s,t)\in \Q\rtimes_\alpha \Q_+^*$,
\[
(x,y)^{-1}(s,t)=\left(\frac{s-x}{y},\frac{t}{y}\right).
\]
It follows that if $(x,y)\in \Q\rtimes_\alpha \Q_+^*$ and $(s,t)\in \N\rtimes_\alpha \N_+^*$, then 
\[
(x,y)^{-1}(s,t)\in \N\rtimes_\alpha \N_+^* \Leftrightarrow t\in y\N_+ \text{ and } s\in x+y\N.
\]
Consequently, the semidirect product $(\Q\rtimes_\alpha \Q_+^*, \N\rtimes_\alpha \N_+^*)$ is a quasi-lattice ordered group with 
\begin{align*}
(x,y)&\vee (u,v)\\
&=\begin{cases}
\big(\min\left((x+y\mathbb{N})\cap(u+v\mathbb{N})\right),\min\left(y\mathbb{N}_+\cap v\mathbb{N}_+\right)\big) & \text{if }(x+y\mathbb{N})\cap(u+v\mathbb{N})\neq\emptyset \\
\infty & \text{otherwise}
\end{cases}
\end{align*}
for any $(x,y), (u,v)\in \Q\rtimes_\alpha \Q_+^*$. It is easy to see that the subsemigroup $\N\rtimes_\alpha \N_+^*$ is not directed: for example $(0,2)\vee (1,2)=\infty$. Let $\mathcal{G}:=\Q\rtimes_\alpha \Q_+^*$ with the discrete topology and let $c:\mathcal{G}\rightarrow \Q\rtimes_\alpha \Q_+^*$ be the identity map. For each $(m,n)\in \N\rtimes_\alpha \N_+^*$, we have that $\mathbf{X}(\mathcal{G})_{(m,n)}=\overline{C_c(\{(m,n)\})}\cong \C$. Thus, each $\phi_{(m,n)}$ is injective and takes values in $\mathcal{K}_{C_r^*(\mathcal{G}_e)}(\mathbf{X}(\mathcal{G})_{(m,n)})$. Hence, $\mathbf{X}(\mathcal{G})$ is compactly aligned by \cite[Proposition~5.8]{MR1907896}. Since $r(\mathcal{G}_{(m,n)})=\{(0,1)\}$ for each $(m,n)\in \N\rtimes_\alpha \N_+^*$, the groupoid $\mathcal{G}$ does not satisfy \eqref{extendibility}. Hence, by Lemma~\ref{sufficient conditions for Nica covariance}, the inclusion of $\mathbf{X}(\mathcal{G})$ in $C_r^*(\mathcal{G})$ is not Nica covariant. We now show that the inclusion map is also not Cuntz--Pimsner covariant. Let $F:=\{(0,1),(0,2),(1,2)\}\subseteq \N\rtimes_\alpha \N_+^*$ and define $T_{(0,1)}:=\phi_{(0,1)}(1)$, $T_{(0,2)}:=\phi_{(0,2)}(-1)$, $T_{(1,2)}:=\phi_{(1,2)}(-1)$. Then $I^{(0,1)}(T_{(0,1)})+I^{(0,2)}(T_{(0,2)})+I^{(1,2)}(T_{(1,2)})=-1\neq 0$. However, as we shall now show, $\iota_{(0,1)}^p(T_{(0,1)})+\iota_{(0,2)}^p(T_{(0,2)})+\iota_{(1,2)}^p(T_{(1,2)})=0$ for large $p$. Let $(x,y)\in \N\rtimes_\alpha \N_+^*$ and choose $z\in 2y\N_+$ and $w\in x+y\N$. Thus, $(w,z)\geq (x,y)$. If $(u,v)\geq (w,z)$, then $v\in z\N_+\subseteq 2\N_+$, and so either $(u,v)\geq (0,2)$ (if $u\in 2\N$) or $(u,v)\geq (1,2)$ (if $u\in 1+2\N$). Since $(0,2)\vee (1,2)=\infty$, we have that
\begin{align*}
\iota_{(0,1)}^{(u,v)}(T_{(0,1)})+\iota_{(0,2)}^{(u,v)}(T_{(0,2)})+\iota_{(1,2)}^{(u,v)}(T_{(1,2)})
&=
\begin{cases}
\iota_{(0,1)}^{(u,v)}(T_{(0,1)})+\iota_{(0,2)}^{(u,v)}(T_{(0,2)}) & \text{if $u\in 2\N$}\\
\iota_{(0,1)}^{(u,v)}(T_{(0,1)})+\iota_{(1,2)}^{(u,v)}(T_{(1,2)}) & \text{if $u\in 1+2\N$}
\end{cases}\\
&=\phi_{(u,v)}(1)+\phi_{(u,v)}(-1)\\
&=0.
\end{align*}
\end{exam}

Combining Lemma~\ref{sufficient conditions for Nica covariance} and Theorem~\ref{induced homomorphism is injective - covariance algebra version}, we get the following. 

\begin{thm}
\label{main result for CNP algebras}
Suppose that $(G,P)$ is a quasi-lattice ordered group and $G$ is amenable. Additionally, suppose that $\mathbf{X}(\mathcal{G})$ is $\tilde{\phi}$-injective, and that either the left actions on the fibres of $\mathbf{X}(\mathcal{G})$ are all injective or that $P$ is directed. Then $\mathbf{X}(\mathcal{G})$ is compactly aligned and the inclusion of $\mathbf{X}(\mathcal{G})$ in $C_r^*(\mathcal{G})$ induces an isomorphism from $\mathcal{NO}_{\mathbf{X}(\mathcal{G})}$ to $C_r^*(\mathcal{G})$ if and only if $\mathcal{G}$ satisfies \eqref{extendibility} and $c^{-1}(P)$ generates $\mathcal{G}$ as a groupoid.
\begin{proof}
Recall that if $\mathbf{X}(\mathcal{G})$ is compactly aligned, then \cite[Proposition~4.6]{SEHNEM2018} tells us that $\mathcal{NO}_{\mathbf{X}(\mathcal{G})}$ and $C_r^*(\mathcal{G}_e)\times_{\mathbf{X}(\mathcal{G})}P$ are canonically isomorphic.  Observe that if $\mathcal{G}$ satisfies \eqref{extendibility}, then Lemma~\ref{without closure, things are simpler} tells us that
\[
r(\mathcal{G}_p)\cap r(\mathcal{G}_{p_1}\mathcal{G}_{p_2}^{-1}\cdots \mathcal{G}_{p_{n-1}}^{-1}\mathcal{G}_{p_n})\subseteq \bigcup_{m\in pP\cap p_1p_2^{-1}\cdots p_{n-1}^{-1}p_nP}r(\mathcal{G}_m) \quad \text{for $p,p_1,\ldots,p_n\in P$},
\]
and so $\mathcal{G}$ satisfies \eqref{necessary and sufficient condition equivalent}. Hence, it follows from Theorem~\ref{induced homomorphism is injective - covariance algebra version} and Lemma~\ref{sufficient conditions for Nica covariance} that $\mathbf{X}(\mathcal{G})$ is compactly aligned and the inclusion of $\mathbf{X}(\mathcal{G})$ in $C_r^*(\mathcal{G})$ induces an isomorphism from $\mathcal{NO}_{\mathbf{X}(\mathcal{G})}$ to $C_r^*(\mathcal{G})$ if and only if $\mathcal{G}$ satisfies \eqref{extendibility} and $c^{-1}(P)$ generates $\mathcal{G}$ as a groupoid.
\end{proof}
\end{thm}

\begin{rem}
We point out that the technical condition of $\mathbf{X}(\mathcal{G})$ being $\tilde{\phi}$-injective is automatic when $(G,P)=(\mathbb{Z}^k,\mathbb{N}^k)$ since $(\mathbb{Z}^k,\mathbb{N}^k)$ satisfies \eqref{maximal element condition}. Moreover, $\mathbb{N}^k$ is directed and $\mathbb{Z}^k$ is amenable, so Theorem~\ref{main result for CNP algebras} is applicable to groupoids with unperforated $\mathbb{Z}^k$-gradings.
\end{rem}

\begin{rem}
We point out that if $\mathcal{G}$ satisfies \eqref{extendibility}, then $c^{-1}(P)$ generating $\mathcal{G}$ as a groupoid is equivalent to
\begin{equation}
\label{condition for surjectivity}
\text{for all $\gamma\in \mathcal{G}$ there exist $\mu,\nu\in c^{-1}(P)$ such that $\gamma=\mu\nu^{-1}$.}
\end{equation}
Clearly, if $\mathcal{G}$ satisfies \eqref{condition for surjectivity}, then $c^{-1}(P)$ generates $\mathcal{G}$ as a groupoid, so we just need to worry about the converse. Suppose that $c^{-1}(P)$ generates $\mathcal{G}$ as a groupoid and $\mathcal{G}$ satisfies \eqref{extendibility}. Since $P$ is a semigroup, in order to show that \eqref{condition for surjectivity} holds, we need only show that if $\mu,\nu \in c^{-1}(P)$, then $\nu^{-1}\mu=\rho \tau^{-1}$ for some $\rho,\tau\in c^{-1}(P)$. Suppose that $\mu\in\mathcal{G}_p$ and $\nu\in\mathcal{G}_q$ for some $p,q\in P$ and $\nu^{-1}\mu$ is defined. Thus, $r(\mu)=r(\nu)\in r(\mathcal{G}_p)\cap r(\mathcal{G}_q)$. By \eqref{extendibility} we have that $p\vee q<\infty$ and there exists $\sigma\in \mathcal{G}_{p\vee q}$ with $r(\mu)=r(\nu)=r(\sigma)$. Then
\[
\nu^{-1}\mu=\nu^{-1}r(\sigma)\mu=\nu^{-1}\sigma \sigma^{-1}\mu=(\nu^{-1}\sigma)(\mu^{-1}\sigma)^{-1}
\in \mathcal{G}_{q^{-1}(p\vee q)}\mathcal{G}_{p^{-1}(p\vee q)}^{-1}
\]
as required.

However, in general $c^{-1}(P)$ generating $\mathcal{G}$ need not imply \eqref{condition for surjectivity}. Recall Example~\ref{when things can go wrong}. We first show that $c^{-1}(\N\rtimes_\alpha \N_+^*)=\N\rtimes_\alpha \N_+^*$ generates $\Q\rtimes_\alpha \Q_+^*$ as a group(oid). Let $a,b\in \mathbb{Z}$ (with $b\neq 0$) and $c,d\in \mathbb{N}_+$ so that $\left(\frac{a}{b},\frac{c}{d}\right)\in \Q\rtimes_\alpha \Q_+^*$. We claim that $\left(\frac{a}{b},\frac{c}{d}\right)=(x,y)^{-1}(w,z)=\big(\frac{w-x}{y},\frac{z}{y}\big)$ for some $(x,y),(w,z)\in \N\rtimes_\alpha \N_+^*$. Without loss of generality, we may assume that $b>0$ (if $b<0$, then just replace $b$ by $-b$ and $a$ by $-a$). Let $y:=bd\in \mathbb{N}_+$ and $z:=bc\in \mathbb{N}_+$. Then $\frac{z}{y}=\frac{c}{d}$. If $a\geq 0$, let $x:=0$ and $w:=ad\in \mathbb{N}$, so that $\frac{w-x}{y}=\frac{ad}{bd}=\frac{a}{b}$. On the other hand, if $a<0$, let $w:=0$ and $x:=-ad\in \mathbb{N}$, so that $\frac{w-x}{y}=\frac{-(-ad)}{bd}=\frac{a}{b}$. Thus, $c^{-1}(\N\rtimes_\alpha \N_+^*)$ generates $\Q\rtimes_\alpha \Q_+^*$ as a group(oid). However, as we now show, this example does not satisfy \eqref{condition for surjectivity}. To see this suppose that there exist $(x,y),(w,z)\in \N\rtimes_\alpha \N_+^*$ such that $\left(\frac{1}{3},\frac{1}{2}\right)=(x,y)(w,z)^{-1}=\big(x-\frac{wy}{z},\frac{y}{z}\big)$. Then $z=2y$, and so we must have that $2=6x-3w$. But this is impossible because the linear Diophantine equation $6x-3w=2$ cannot have an integer solution because $\mathrm{gcd}(6,3)=3$ does not divide $2$.
\end{rem}


\section{Examples and Applications}
\label{examples}

In addition to our standing hypotheses from the start of Section~\ref{constructing the product system}, we now have a number of additional conditions that our groupoid and cocycle need to satisfy in order to apply Theorem~\ref{induced homomorphism is injective - covariance algebra version} or Theorem~\ref{main result for CNP algebras}. We now present a couple of examples that our theorems can be applied to. 

\subsection{Topological higher-rank graphs}

We summarise the relevant background material on topological higher-rank graphs and their associated groupoids from \cite{MR2301938}. A topological $k$-graph consists of a small category $\Lambda$ and a functor $d:\Lambda\rightarrow \mathbb{N}^k$ (called the degree functor), such that
\begin{enumerate}[label=\upshape(\roman*)]
\item
$\mathrm{Obj}(\Lambda)$, $\mathrm{Mor}(\Lambda)$ are second-countable locally compact Hausdorff spaces;
\item
$r,s:\mathrm{Mor}(\Lambda)\rightarrow \mathrm{Obj}(\Lambda)$ are continuous and $s$ is a local homeomorphism;
\item
the composition map from $\Lambda\times_c\Lambda:=\{(\lambda,\mu)\in \Lambda\times \Lambda: s(\lambda)=r(\mu)\}$ to $\Lambda$ is continuous and open;
\item
$d$ is continuous (equipping $\mathbb{N}^k$ with the discrete topology);
\item
for all $\lambda\in \Lambda$ and $m,n\in \mathbb{N}^k$ such that $d(\lambda)=m+n$, there exists a unique $(\mu,\nu)\in \Lambda\times_c\Lambda$ with $d(\mu)=m$ and $d(\nu)=n$ such that that $\lambda=\mu\nu$.
\end{enumerate}
It follows from (v) that if $\lambda\in\Lambda$ and $m,n\in \mathbb{N}^k$ with $0\leq m\leq n\leq d(\lambda)$, then there exist unique $\mu,\nu,\rho\in \Lambda$ with $d(\mu)=m$, $d(\nu)=n-m$, and $d(\rho)=d(\lambda)-n$ such that $\lambda=\mu\nu\rho$. We write $\lambda(0,m)$, $\lambda(m,n)$, and $\lambda(n,d(\lambda))$ for $\mu$, $\nu$, and $\rho$ respectively.

For $n\in \mathbb{N}^k$, write $\Lambda^n:=d^{-1}(n)$. We define $\Lambda \times_s \Lambda:=\{(\lambda,\mu)\in \Lambda\times \Lambda:s(\lambda)=s(\mu)\}$. Given subsets $U,V\subseteq \Lambda$, we write $UV:=\{\lambda\mu:\lambda\in U, \mu\in V, s(\lambda)=r(\mu)\}$ and $U\times_sV:=(U\times V) \cap(\Lambda \times_s \Lambda)$.

We say that $\Lambda$ is compactly aligned if for all $m,n\in \mathbb{N}^k$ and compact sets $U\subseteq \Lambda^m$ and $V\subseteq \Lambda^n$, the set
\[
U\vee V:=U\Lambda^{m\vee n-m}\cap V\Lambda^{m\vee n-n}
\]
is compact.

We now describe how to associate a groupoid to a topological $k$-graph $\Lambda$. We begin by introducing the path space of $\Lambda$. For $k\in \mathbb{N}$ and $m\in (\mathbb{N}\cup\{\infty\})^k$, we define the topological $k$-graph $\Omega_{k,m}$ by equipping
\[
\mathrm{Obj}(\Omega_{k,m}):=\{p\in \mathbb{N}^k:p\leq m\}
\]
and
\[
\mathrm{Mor}(\Omega_{k,m}):=\{(p,q)\in \mathbb{N}^k\times \mathbb{N}^k:p\leq q\leq m\}
\]
with the discrete topologies and defining
\[
r(p,q):=p, \quad s(p,q):=q, \quad (p,q)(q,r):=(p,r), \quad d(p,q):=q-p.
\]
We say that a continuous functor from $\Omega_{k,m}$ to $\Lambda$ is a graph morphism if it preserves the degree functor. The path space of $\Lambda$ is then
\[
X_\Lambda:=\bigcup_{m\in (\mathbb{N}\cup\{\infty\})^k}\{x:\Omega_{k,m}\rightarrow \Lambda \text{ is a graph morphism}\}.
\]
Given $\lambda\in \Lambda$, there exists a unique graph morphism $x_\lambda:\Omega_{k,d(\lambda)}\rightarrow \Lambda$ such that $x_\lambda(0,d(\lambda))=\lambda$. Consequently, we may view $\Lambda$ as a subset of $X_\Lambda$. We also extend the range and degree maps to $X_\Lambda$ by setting $r(x):=x(0)$ and $d(x):=m$ for a graph morphism $x:\Omega_{k,m}\rightarrow \Lambda$.

We can add/remove finite paths to/from the end of elements of $X_\Lambda$ as follows. If $x\in X_\Lambda$ and $\lambda\in \Lambda$ with $s(\lambda)=r(x)$, then there exists a unique $\lambda x\in X_\Lambda$ with $d(\lambda x)=d(\lambda)+d(x)$ such that
\[
(\lambda x)(0,p)=
\begin{cases}
\lambda(0,p) & \text{if $p\leq d(\lambda)$}\\
\lambda x(0,p-d(\lambda)) & \text{if $d(\lambda)\leq p\leq d(\lambda x)$.}
\end{cases}
\]
If $x\in X_\Lambda$ and $m\in \mathbb{N}^k$ with $m\leq d(x)$, then there exists a unique $\sigma^m x\in X_\Lambda$ with $d(\sigma^m x)=d(x)-m$ such that
\[
(\sigma^m x)(0,p)=x(m,m+p) \quad \text{for $p\leq d(\sigma^m x)$.}
\]

Finally, we are ready to associate a groupoid to $\Lambda$. We set
\begin{align*}
G_\Lambda:
&=
\{(\lambda x,d(\lambda)-d(\mu),\mu x)\in X_\Lambda\times \Z^k\times X_\Lambda: \lambda,\mu\in \Lambda, x\in X_\Lambda, s(\lambda)=s(\mu)=r(x)\}
\\
&=
\left\{ (x,m,y)\in X_\Lambda\times \Z^k\times X_\Lambda : \begin{array}{l}
    \text{there exists } p,q\in \N^k \text{ such that } m=p-q, \\
     d(x)\geq p, d(y)\geq q, \text{ and }\sigma^p x=\sigma^q y
  \end{array}\right\}.
\end{align*}
Then $G_\Lambda$ has the structure of a groupoid, which we call the path groupoid, with operations
\[
(x,m,y)(y,n,z):=(x,m+n,z) \quad \text{and} \quad (x,m,y)^{-1}:=(y,-m,x).
\]
We identify the unit space $G_\Lambda^{(0)}$ with $X_\Lambda$ via the map $(x,0,x)\mapsto x$.

We define a topology on $G_\Lambda$ as follows. For $m\in \mathbb{N}^k$ and $F\subseteq \Lambda\times_s \Lambda$, we define $Z(F,m)\subseteq G_\Lambda$ by
\[
Z(F,m):=\{(\lambda x,d(\lambda)-d(\mu),\mu x)\in G_\Lambda:(\lambda,\mu)\in F, d(\lambda)-d(\mu)=m\}.
\]
The collection of sets $Z(U\times_s V,m)\setminus Z(F,m)$, where $m\in \mathbb{N}^k$, $U,V\subseteq \Lambda$ are open, and $F\subseteq \Lambda\times_s \Lambda$ is compact, forms a basis for a second-countable Hausdorff topology on $G_\Lambda$. If $\Lambda$ is compactly aligned, then this topology is locally compact, and gives $G_\Lambda$ the structure of a second-countable locally compact Hausdorff \'{e}tale groupoid.

We note that there exists a cocycle $c:G_\Lambda\rightarrow \Z^k$ given by $c(x,m,y):=m$. It is not hard to see that $c$ is continuous --- for any $m\in \mathbb{N}^k$, $c^{-1}(m)=Z(\Lambda\times_s \Lambda,m)\setminus Z(\emptyset,m)$, which is a basic open set.

Given a compactly aligned topological $k$-graph, we also define another groupoid (called the boundary-path groupoid) as a certain reduction of the path groupoid $G_\Lambda$. First, we need to introduce the notion of an exhaustive set.

For $V\subseteq \Lambda^0$, we say that $E\subseteq V\Lambda$ is exhaustive for $V$ if for all $\lambda\in V\Lambda$ there exists $\mu\in E$ such that $\{\lambda\}\Lambda^{d(\lambda)\vee d(\mu)-d(\lambda)}\cap \{\mu\}\Lambda^{d(\lambda)\vee d(\mu)-d(\mu)}\neq \emptyset$. For $v\in \Lambda^0$, let $v\mathrm{CE}(\Lambda)$ denote the set of all compact sets $E\subseteq \Lambda$ such that $r(E)$ is a neighbourhood of $v$ and $E$ is exhaustive for $r(E)$.

We say that a path $x\in X_\Lambda$ is a boundary-path if for all $m\in \mathbb{N}^k$ with $m\leq d(x)$, and for all $E\in x(m)\mathrm{CE}(\Lambda)$, there exists $\lambda\in E$ such that $x(m,m+d(\lambda))=\lambda$. We write $\partial\Lambda$ for the set of all boundary paths. One can then show that $\partial\Lambda$ is a closed invariant subset of $G_\Lambda^{(0)}=X_\Lambda$ (i.e. if $x\in \partial \Lambda$, and $m\in \mathbb{N}^k$ with $m\leq d(x)$ and $\lambda\in \Lambda$ with $s(\lambda)=r(x)$, then $\sigma^m(x), \lambda x\in \partial \Lambda$). Consequently, we can define the reduction $\mathcal{G}_\Lambda:=G_\Lambda|_{\partial\Lambda}$. We call $\mathcal{G}_\Lambda$ the boundary-path groupoid. The cocycle $c:G_\Lambda\rightarrow \Z^k$ descends to a continuous $\mathbb{Z}^k$-valued cocycle on $\mathcal{G}_\Lambda$, which we denote by $c'$.

\begin{exam}
\label{topological higher rank graphs}
Let $\Lambda$ be a compactly aligned topological $k$-graph, and $G_\Lambda$ and $\mathcal{G}_\Lambda$ be the associated path and boundary-path groupoids. Then the cocycle $c:G_\Lambda\rightarrow \mathbb{Z}^k$ is unperforated, $G_\Lambda$ satisfies \eqref{extendibility}, and $c^{-1}(\mathbb{N}^k)$ generates $G_\Lambda$ as a groupoid. Similarly, the cocycle $c':\mathcal{G}_\Lambda\rightarrow \mathbb{Z}^k$ is unperforated, $\mathcal{G}_\Lambda$ satisfies \eqref{extendibility}, and $c'^{-1}(\mathbb{N}^k)$ generates $\mathcal{G}_\Lambda$. Consequently, by Theorem~\ref{main result for CNP algebras} both $C_r^*(G_\Lambda)$ and $C_r^*(\mathcal{G}_\Lambda)$ may be realised as the Cuntz--Nica--Pimsner algebras of compactly aligned product systems over $\N^k$ with coefficient algebras $C_r^*((G_\Lambda)_0)$ and $C_r^*((\mathcal{G}_\Lambda)_0)$ respectively.
\begin{proof}
Firstly, let us consider the path groupoid $G_\Lambda$. We begin by showing that $c$ is unperforated. Suppose $\gamma\in (G_\Lambda)_{m+n}$ for some $m,n\in \N^k$. Thus, $\gamma=(x,m+n,y)$ for some $x,y\in X_\Lambda$, and there exist $p,q\in \N^k$ such that $m+n=p-q$, $d(x)\geq p$, $d(y)\geq q$, and $\sigma^p(x)=\sigma^q(y)$. Since $d(x)\geq p=m+n+q\geq m$, we see that $\sigma^m(x)$ is well-defined. Furthermore, since $m=p-(q+n)$ and $\sigma^p(x)=\sigma^{q+n}(\sigma^m(x))$, we see that $(x,m,\sigma^m(x))\in (G_\Lambda)_m$. Also, since $n=(p-m)-q$ and $\sigma^{p-m}(\sigma^m(x))=\sigma^p(x)=\sigma^q(y)$, we see that $(\sigma^m(x),n,y)\in (G_\Lambda)_n$. As $\gamma=(x,m+n,y)=(x,m,\sigma^m(x))(\sigma^m(x),n,y)$, we conclude that $c$ is unperforated.

It is not difficult to see that $G_\Lambda$ satisfies \eqref{extendibility}. Let $m,n\in \N^k$. Since $c$ is unperforated, we know that $r((G_\Lambda)_{m\vee n})\subseteq r((G_\Lambda)_m)\cap r((G_\Lambda)_n)$, and so we just need to establish the reverse containment. Suppose that $x\in r((G_\Lambda)_m)\cap r((G_\Lambda)_n)$. Hence, $x\in X_\Lambda$ and there exist $y,z\in X_\Lambda$ such that $(x,m,y)\in (G_\Lambda)_m$ and $(x,n,z)\in (G_\Lambda)_n$. Thus, there exist $p,q,s,t\in \N^k$ such that $m=p-q$, $n=s-t$, and $d(x)\geq p,s$. Since $d(x)\geq p=m+q\geq m$ and $d(x)\geq s=n+t\geq n$, we have that $d(x)\geq m\vee n$. Hence, $\sigma^{m\vee n}(x)$ is well-defined and $(x,m\vee n, \sigma^{m\vee n}(x))\in (G_\Lambda)_{m\vee n}$. Thus, $x\in r((G_\Lambda)_{m\vee n})$, and we conclude that \eqref{extendibility} holds.

We now show that $G_\Lambda$ is generated by $c^{-1}(\mathbb{N}^k)$. Let $(x,m,y)\in G_\Lambda$. Hence there exist $p,q\in \N^k$ such that $m=p-q$, $d(x)\geq p$, $d(y)\geq q$, and $\sigma^p(x)=\sigma^q(y)$. Thus, $(x,p,\sigma^p(x))\in (G_\Lambda)_p$ and $(y,q,\sigma^q(y))\in (G_\Lambda)_q$. Since $(x,m,y)=(x,p-q,y)=(x,p,\sigma^p(x))(y,q,\sigma^q(y))^{-1}$, we conclude that $c^{-1}(\mathbb{N}^k)$ generates $G_\Lambda$.

It remains to consider the boundary-path groupoid. Recall that if $x\in \partial \Lambda$ and $m\leq d(x)$, then $\sigma^m(x)\in \partial \Lambda$. Hence, the same working as for the path groupoid $G_\Lambda$ shows that the cocycle $c':\mathcal{G}_\Lambda\rightarrow \Z^k$ is unperforated, $\mathcal{G}_\Lambda$ satisfies \eqref{extendibility}, and $c'^{-1}(\mathbb{N}^k)$ generates $\mathcal{G}_\Lambda$.
\end{proof}
\end{exam}

\subsection{Semigroup action groupoids}

We can generalise the previous example by considering semigroup action groupoids \cite[\S5]{MR3592511}. For those interested in the relationship between semigroup action groupoids and topological higher-rank graphs, see \cite[\S6]{MR3592511} (in fact the more general situation of topological $P$-graphs is considered). We now briefly summarise the necessary background and definitions for semigroup action groupoids.

Let $X$ be a set and $P$ a unital semigroup. A right partial action of $P$ on $X$ consists of a subset $X*P\subseteq X\times P$ and a map $T:X*P\rightarrow X$ that sends $(x,p)$ to $x\cdot p$, satisfying the following conditions
\begin{enumerate}[label=\upshape(\roman*)]
\item
for all $x\in X$, $(x,e)\in X*P$ and $x\cdot e=x$;
\item
for all $(x,p,q)\in X\times P\times P$, $(x,pq)\in X*P$ if and only if $(x,p)\in X*P$ and $(x\cdot p,q)\in X*P$; if this holds then $x\cdot (pq)=(x\cdot p)\cdot q$.
\end{enumerate}

For $p\in P$, write $U(p):=\{x\in X:(x,p)\in X*P\}$, $V(p):=\{x\cdot p\in X:(x,p)\in X*P\}$, and define $T_p:U(p)\rightarrow V(p)$ by $T_p(x):=x\cdot p$. It follows from condition (ii) that if $p,q\in P$, then $U(pq)\subseteq U(p)$. We call the triple $(X,P,T)$ a semigroup action. We say that $(X,P,T)$ is directed if for all $p,q\in P$ such that $U(p)\cap U(q)\neq \emptyset$, there exists $r\in pP\cap qP$ such that $U(p)\cap U(q)=U(r)$.

As discussed in \cite[Example~5.3]{MR3592511}, if $(X,P,T)$ is a semigroup action and $X*P=X\times P$ (i.e. the action is everywhere defined), then the semigroup action is directed if and only if $P$ is directed (in the sense that $pP\cap qP\neq \emptyset$ for each $p,q\in P$).
If $P$ is a subsemigroup of a group $G$, then $P$ is directed if and only if $P$ satisfies the Ore condition $P^{-1}P\subseteq PP^{-1}$.

Now assume that $P$ is a subsemigroup of a group $G$ and $(X,P,T)$ is a directed action. We let $\mathcal{G}(X,P,T)$ denote the collection of triples
\[
\{(x,g,y)\in X\times G\times X:\exists m,n\in P \text{ s.t. } g=mn^{-1},\, x\in U(m), y\in U(n),\, x\cdot m=y\cdot n\}.
\]
Then $\mathcal{G}(X,P,T)$ has the structure of a groupoid, which we call the semidirect product groupoid, with operations
\[
(x,g,y)(y,h,z)=(x,gh,z) \quad \text{and} \quad (x,g,y)^{-1}=(y,g^{-1},x)
\]
for $x,y,z\in X$, $g,h\in G$. The semidirect product groupoid carries a canonical cocycle $c:\mathcal{G}(X,P,T)\rightarrow G$ defined by $c(x,g,y):=g$.

We say that a semigroup action $(X,P,T)$ is locally compact if
\begin{enumerate}[label=\upshape(\roman*)]
\item
$X$ is a locally compact Hausdorff space;
\item
$G$ is a discrete group;
\item
for each $p\in P$, $U(p)$ and $V(p)$ are open subsets of $X$ and $T_p:U(P)\rightarrow V(P)$ is a local homeomorphism.
\end{enumerate}

Given a locally compact semigroup action $(X,P,T)$ we define a topology on $\mathcal{G}(X,P,T)$ as follows. For $p,q\in P$ and sets $U,V\subseteq X$, we let
\[
\mathcal{Z}(U,p,q,V):=\{(x,pq^{-1},y)\in \mathcal{G}(X,P,T):x\in U, y\in V, x\cdot p=y\cdot q\}.
\]
The collection of sets $\mathcal{Z}(U,p,q,V)$, where $p,q\in P$ and $U,V$ are open subsets of $X$ forms a basis for a topology on $\mathcal{G}(X,P,T)$. This topology gives $\mathcal{G}(X,P,T)$ the structure of a locally compact Hausdorff \'{e}tale groupoid. With respect to this topology the cocycle $c:\mathcal{G}(X,P,T)\rightarrow G$ is continuous. If $X$ is second-countable and $P$ is countable, then this topology on $\mathcal{G}(X,P,T)$ is also second-countable.

\begin{exam}
\label{semigroup action groupoids}
Let $G$ be a group and $P\subseteq G$ a unital subsemigroup. Suppose that $(X,P,T)$ is a locally compact directed semigroup action and let $\mathcal{G}(X,P,T)$ be the associated semidirect product groupoid. Then the continuous cocycle $c:\mathcal{G}(X,P,T)\rightarrow G$ defined by $c(x,g,y):=g$ is unperforated and $c^{-1}(P)$ generates $\mathcal{G}(X,P,T)$ as a groupoid. Furthermore,
\begin{equation}
\label{range condition}
r(\mathcal{G}(X,P,T)_p)\cap r(\mathcal{G}(X,P,T)_q)= \bigcup_{r\in pP\cap qP}r(\mathcal{G}(X,P,T)_r)
\quad \text{for any $p,q\in P$}.
\end{equation}
Hence, by Lemma~\ref{without closure, things are simpler}, $\mathcal{G}(X,P,T)$ satisfies \eqref{necessary and sufficient condition equivalent} (and if $(G,P)$ is quasi-lattice ordered, $\mathcal{G}(X,P,T)$ satisfies \eqref{extendibility}). Thus, if $G$ is an amenable group, Theorem~\ref{induced homomorphism is injective - covariance algebra version} tells us that $C_r^*(\mathcal{G}(X,P,T))$ may be realised as the covariance algebra of a product system over $P$ with coefficient algebra $C_r^*(\mathcal{G}(X,P,T)_e)$. Moreover, if $G$ is an amenable group and $(G,P)$ is quasi-lattice ordered, then Theorem~\ref{main result for CNP algebras} tells us that $C_r^*(\mathcal{G}(X,P,T))$ may be realised as the Cuntz--Nica--Pimsner algebra of a compactly aligned product system over $P$ with coefficient algebra $C_r^*(\mathcal{G}(X,P,T)_e)$

\begin{proof}
We begin by showing that the cocycle $c$ is unperforated. Suppose $(x,mn,y)\in \mathcal{G}(X,P,T)$ where $m,n\in P$. Hence there exist $s,t\in P$ such that $mn=st^{-1}$ and $x\in U(s)$, $y\in U(t)$ with $x\cdot s=y\cdot t$. Since $nt\in P$ and $s=m(nt)$, we have that $x\in U(s)\subseteq U(m)$, and so $(x,m,x\cdot m)\in \mathcal{G}(X,P,T)_m$. Furthermore, as $y\cdot t=x\cdot s=(x\cdot m)\cdot (nt)$, we see that $x\cdot m\in U(nt)$. Since $n=(nt)t^{-1}$, we see that $(x\cdot m,n,y)\in \mathcal{G}(X,P,T)_n$. Finally since, $(x,mn,y)=(x,m,x\cdot m)(x\cdot m,n,y)$ we conclude that $c$ is unperforated.

Secondly, we check that $c^{-1}(P)$ generates $\mathcal{G}(X,P,T)$ as a groupoid. Let $(x,g,y)\in \mathcal{G}(X,P,T)$. Choose $p,q\in P$ such that $g=pq^{-1}$, $x\in U(p)$, $y\in U(q)$, and $x\cdot p=y\cdot q$. Then $(x,p,x\cdot p)\in \mathcal{G}(X,P,T)_p$ and $(y,q,y\cdot q)\in \mathcal{G}(X,P,T)_q$. Since
\[
(x,g,y)=(x,pq^{-1},y)=(x,p,x\cdot p)(y\cdot q,q^{-1},y)=(x,p,x\cdot p)(y,q,y\cdot q)^{-1},
\]
we conclude that $c^{-1}(P)$ generates $\mathcal{G}(X,P,T)$ as a groupoid.

Finally, we show that $\mathcal{G}(X,P,T)$ satisfies \eqref{range condition}. The inclusion $\bigcup_{r\in pP\cap qP}r(\mathcal{G}(X,P,T)_r)\subseteq r(\mathcal{G}(X,P,T)_p)\cap r(\mathcal{G}(X,P,T)_q)$ for $p,q\in P$ follows from the fact that $c$ is unperforated. To see the reverse containment, fix $p,q\in P$ and suppose that $x\in r(\mathcal{G}(X,P,T)_p)\cap r(\mathcal{G}(X,P,T)_q)$. Hence, there exist $y,z\in X$ such that $(x,p,y),(x,q,z)\in \mathcal{G}(X,P,T)$. Choose $m,n,s,t\in P$ such that $p=mn^{-1}$, $q=st^{-1}$, $x\in U(m)$, $y\in U(n)$, $x\in U(s)$, $z\in U(t)$, and $x\cdot m=y\cdot n$, $x\cdot s=z\cdot t$. Since $x\in U(m)\cap U(s)$ and the semigroup action is directed, there exists $r\in mP\cap sP$ such that $x\in U(r)$. Hence, $(x,r,x\cdot r)\in \mathcal{G}(X,P,T)_r$, and we see that $x\in r(\mathcal{G}(X,P,T)_r)$. Since $r\in mP\cap sP=(pn)P\cap (qt)P\subseteq pP\cap qP$, we conclude that $x\in \bigcup_{r\in pP\cap qP}r(\mathcal{G}(X,P,T)_r)$. Hence, $r(\mathcal{G}(X,P,T)_p)\cap r(\mathcal{G}(X,P,T)_q)\subseteq \bigcup_{r\in pP\cap qP}r(\mathcal{G}(X,P,T)_r)$, and we see that \eqref{range condition} holds.
\end{proof}
\end{exam}

\end{document}